\documentclass{article}
\usepackage{amssymb}
\usepackage{amsmath}
\usepackage{color}

\usepackage[square]{natbib}         %%%%round per parentesi tonde%%%
\newtheorem{lem}{Lemma}[section]
\newtheorem{cor}[lem]{Corollary}
\newtheorem{teo}[lem]{Theorem}
\newtheorem{os}[lem]{Remark}
\newtheorem{defi}[lem]{Definition}
\newtheorem{prop}[lem]{Proposition}
\newtheorem{esem}[lem]{Example}

\newcommand{\qed}{\thinspace\null\nobreak\hfill\hbox{\vbox{\kern-.2pt\hrule
 height.2pt depth.2pt\kern-.2pt\kern-.2pt \hbox to2.5mm{\kern-.2pt\vrule
 width.4pt \kern-.2pt\raise2.5mm\vbox to.2pt{}\lower0pt\vtop
 to.2pt{}\hfil\kern-.2pt \vrule
 width.4pt \kern-.2pt}\kern-.2pt\kern-.2pt\hrule height.2pt depth.2pt
 \kern-.2pt}}\par\medbreak}

\newcommand{\R}{\mathbb{R}}

\newcommand{\Ha}{\mathcal{H}}

\newcommand{\N}{\mathbb{N}}

\date{}

\begin{document}

\title{Sample distribution theory using Coarea Formula}
\author{{ 
L. Negro} \thanks{Dipartimento di Matematica e Fisica ``Ennio De Giorgi'', Universit\`a del Salento, C.P.193, 73100, Lecce, Italy.
email:   
luigi.negro@unisalento.it} }

\maketitle
\begin{abstract}
Let $\left(\Omega,\Sigma,p\right)$ be   a probability measure space   and let $X:\Omega\to\R^k$ be  a   (vector valued) random variable. We suppose that the probability $p_X$ induced by $X$ is absolutely continuous with respect to the Lebesgue measure on $\R^k$ and set   $f_X$ as its density function. Let $\phi:\R^k\to \R^n$ be a $C^1$-map and let us consider the new random variable $Y=\phi(X):\Omega\to\R^n$.
Setting $m:=\max\{\mbox{rank }(J\phi(x)):x\in\R^k\}$, we prove  that the  probability  $p_Y$ induced by  $Y$ has a density function $f_Y$ with respect to the  Hausdorff measure $\Ha^m$  on $\phi(\R^k)$ which  satisfies  
\begin{align*}
	f_Y(y)=
	\int_{\phi^{-1}(y)}f_X(x)\frac{1}{J_m\phi(x)}\,d\Ha^{k-m}(x), &\quad \text{for $\Ha^m$-a.e.}\quad y\in\phi(\R^k).
\end{align*} 
Here $J_m\phi$ is the $m$-dimensional Jacobian of $\phi$. When $J\phi$ has maximum rank we allow the map $\phi$   to be only locally Lipschitz.   We also consider the case of  $X$ having   probability  concentrated on  some $m$-dimensional sub-manifold $E\subseteq\R^k$ and provide, besides, several examples including  algebra of random variables, order statistics, degenerate normal distributions, Chi-squared and "Student's t" distributions.

\par

\end{abstract}

\section{Introduction}
The theory of Sample Distribution  is an important branch of  Statistics and Probability which study the general problem of determining the distribution of functions of random vectors. It   provides a formal framework for  modelling, simulating and making statistical inference. To be more precise, let us fix   a probability measure space $\left(\Omega,\Sigma,p\right)$  and let $X$ be  a   (vector valued) random variable    i.e a $\Sigma$-measurable function from $\Omega$ to $\R^k$, where $k\in\N$; $X$ is usually referred to as the data. Let $Y$ be any measurable function of the data $X$ i.e. $Y$ is a  random variable which satisfies $Y=\phi(X)$ for some measurable function $\phi:\R^k\to\R^n$ with  $n\in\N$; $Y$ is usually called a statistic. 
The problem which we address consists in finding the probability distribution of $Y$ knowing the distribution of $X$.

 Depending on the nature of the data, there are in general different approaches for finding the distribution of the statistic $Y$,  including the distribution function technique,  the moment-generating function technique and the change of variable technique (see e.g. \citep[Section 4.1, page 122]{HoggCraig}). In the last case,  let us suppose for example that  the probability measure induced by $X$ is absolutely continuous with respect to the Lebesgue measure on $\R^k$ and let  $f_X$ be its density function. Then if $k=n$ and $\phi:\R^k\to\R^k$ is a $C^1$-diffeomorphism, then the  change of variable $y=\phi(x)$ yields, for every Borel subset $A\subseteq\R^k$,
\begin{align}\label{int 2}
p\left(\phi(X)\in A\right)=p\left(X\in \phi^{-1}(A)\right)=\int_{\phi^{-1}(A)}f_X(x)\,dx=\int_{A}f_X(\phi^{-1}(y))\frac{1}{|\mbox{det}J_\phi(\phi^{-1}(y))|}\,dy.
\end{align}
The last equation implies that the probability measure induced by  $Y$ is absolutely continuous with respect to the Lebesgue measure on $\R^k$ and   its density function $f_Y$ is  determined by the equation
\begin{align}\label{int 1}
f_Y(y)=f_X(\phi^{-1}(y))\frac{1}{|\mbox{det}J_\phi(\phi^{-1}(y))|}, \quad y\in\R^k.
\end{align}
This change of variables formula is widely used, for example,  in machine learning and is essential for some recent results in density estimation and
	generative modeling like normalizing flows (\citep{Rezende}), NICE (\citep{Dinh2}), or Real NVP (\citep{Dinh1}). However  all uses of this formula in the machine learning literature that we are aware of are constrained by the  bijectivity and the differentiability of the map $\phi$.

In this expository paper we extend formula \eqref{int 1} to the more general case of statistics $Y=\phi(X)$ defined by  locally  Lipschitz functions $\phi:\R^k\to\R^n$. The approach presented is mainly based upon the Coarea Formula proved by Federer in \citep{Federer} which provides, in our setting,  an elegant tool to derive the  distribution of $Y$. This also allows  to get rid of the  differentiability and invertibility assumptions on $\phi$ and to treat the case $k\neq n$: this is  quite useful in many  problems of statistical inference and of machine learning (see e.g. \citep{Civitkovic}).

When $\phi:\R^k\to \R^n$ is a $C^1$-map  and  $m:=\max\{\mbox{rank }(J\phi(x)):x\in\R^k\}$, our result states that the  probability  $p_Y$ induced by  $Y=\phi(X)$ has a density function $f_Y$ with respect to the  Hausdorff measure $\Ha^m$  on $\phi(\R^k)$
which  satisfies  
\begin{align}\label{eq intr}
	f_Y(y)=
	\int_{\phi^{-1}(y)}f_X(x)\frac{1}{J_m\phi(x)}\,d\Ha^{k-m}(x), &\quad \text{for $\Ha^m$-a.e.}\quad y\in\phi(\R^k),
\end{align} 
where $J_m\phi(x)$ is the $m$-dimensional Jacobian of $\phi$ (see Definition \ref{k-Jac}).  When  the Jacobian matrix $J\phi$ has, at any point,  maximum rank we allow the map $\phi$   to be only locally Lipschitz. We also consider the case of  $X$ having   probability  concentrated on  some $m$-dimensional sub-manifold $E\subseteq\R^k$. This case has  many applications in directional and axial statistics, morphometrics, medical
diagnostics, machine vision, image analysis and  molecular biology (see e.g \citep{Bhattacharya} and references there in).

 We do not make any attempt to reach any novelty  but the evidence suggests that this result is not as universally known as it should be. Besides this we give also several examples.

\smallskip
Let us briefly describe the content of the sections. In Section \ref{preliminaries} we  introduce the main definitions and notation that we use throughout the paper.  Section \ref{Section AC}    collects the main results we need  about  the theory of Hausdorff measures and the  Area and Coarea Formulas. In Section \ref{Section Distribution} we  develop the method  proving Formula \eqref{eq intr}:  when $J\phi$ has maximum rank we also allow the map $\phi$   to be locally Lipschitz. In Section \ref{Section Generalization}  we gives  a further  generalization considering   the case of   random variables $X$ having   probability density functions $f_X$ with respect to  the Hausdorff measure $\Ha^m_{\vert E}$ concentrated on  some $m$-dimensional sub-manifold $E\subseteq\R^k$. 
%Finally, in Section \ref{Section apll}, we apply the latter results providing several examples  including some formula for  algebra of random variables the derivation of the probability densities of    algebra of random variables and of order statistics, degenerate normal distributions, Chi-squared and "Student's t" distributions 
Finally, in Section \ref{Section apll}, we provide several examples which show how to apply the latter results in order to characterize the distribution of  algebra of random variables and how  to compute, in an easy way, the probability densities of some classic distributions including order statistics, degenerate normal distributions, Chi-squared and "Student's t" distributions.

\bigskip
\noindent\textbf{Notation.} 
%We adopt standard notation for $L^p$ and Sobolev spaces when $1 \le p<\infty$. For $T>0$ we define $S_T:=(0,T]\times\R^N$. If $a\in\R$ we write $a^+=\max\{a,0\}$ , $a^-=\min\{a,0\}$. 
We write $\langle \lambda,\mu \rangle=\sum_{i}\lambda_i\mu_i$ to denote the inner product of $\R^k$.
 When $f:\R^k\to\R^n$ is a Lipschitz map we write $Jf$ to denote its  Jacobian matrix $\left(\frac{\partial f_i}{\partial x_j}\right)_{i,j}$ which is defined a.e. on $\R^k$. When $A=(a_{ij})\in \R^{n,k}$ is a real matrix we write $Ax$ to denote the linear operator 
$$\phi:\R^k\to\R^n,\qquad x=(x_1,\dots,x_k)\mapsto x\cdot A^t=(y_1,\dots,y_n),\quad y_i=\sum_{j=1}^k a_{i,j}x_j.$$
 With this notation the Jacobian matrix $J(Ax)$ of $Ax$ satisfies $J(Ax)=A$. $I_k$ is the identity matrix of $\R^{k,k}$. If  $\left(\Omega_1, \Sigma_1\right)$ and $\left(\Omega_2, \Sigma_2\right)$ are measurable spaces, a function $f :\Omega_1\to\Omega_2$ is said to be $\left(\Sigma_1,\Sigma_2\right)$-measurable if $f^{-1}(B)\in\Sigma_1$ for all $B\in\Sigma_2$. 
Unless otherwise specified when $\Omega_2=\R^k$ we always choose $\Sigma_2$ as the $\sigma$-algebra $\mathcal{B}(\R^k)$ of all the Borel subsets of $\R^k$ and in this case we simply say that $f$ is  $\Sigma_1$-measurable. We finally write $\mathcal L^k$ and $\mathcal H^s$ to denote respectively the Lebesgue measure and the $s$-dimensional Hausdorff measure on $\R^k$: under this notation we have in particular that  $\mathcal L^k=\mathcal H^k$ and that $\mathcal H^0$ is the counting measure on $\R^k$. 
%$\sigma_{i,j}=\mbox{Cov}\left(X_i,X_j\right)$
%
%$E(X)=$

% \textcolor{red}{a.e. means a.e. respect Lebesgue, $L^1(\R^k)$,  For a subset $F$ we write $F^c$ to denote its complementary set.}\\
% If $M$ has real eigenvalues, $\lambda_{max}(M)$ and $\lambda_{min}(M)$ will denote respectively the maximum and the minimum eigenvalue of $M$.
%
%\bigskip
%\noindent\textbf{Acknowledgement.} The author is grateful to Prof. G. Metafune for several discussions on the topic and
%for his useful comments on this paper.
\bigskip

\noindent\textbf{Acknowledgements.} The author thanks F. Durante  for several comments on a previous version of the manuscript.

\section{Preliminaries}\label{preliminaries}
In this section we fix the main notation and collect the main  results we use concerning the  Probability theory. For a good survey on the topic  we refer  the reader, for example,  to \citep[Chapter IV]{HalmosMeasure} and \citep[Appendix A and B]{Schervish}.

Let $\mu,\nu$ two (positive) measures defined on a  measurable space $\left(\Omega,\Sigma\right)$.   $\nu$ is said to be \emph{absolutely continuous} with respect to $\mu$ and we write $\nu\ll\mu$ if and only if $\nu(B)=0$ for every $B\in\Sigma$ such that $\mu(B)=0$.  $\nu$ is said to have a  \emph{density function} $f$ with respect to $\mu$  if and only if  there exists a measurable positive function  $f:\Omega\to \R^+$ such that 
\begin{align*}
\nu(B)=\int_B f\,d\mu,\quad  \text{for all}\quad A\in\Sigma,
\end{align*}
(note that $f$ is uniquely defined up to zero measure   sets). When $\mu$ is $\sigma$-finite, thanks to the Radon-Nikodym Theorem,  the latter two definitions coincide and  $\frac{d\nu}{d \mu}:=f$ is called the \emph{Radon-Nikodym derivative} of $\nu$ with respect to $\mu$ (see e.g. \citep[Theorem 1.28]{AFPallara}).

Let now $\left(\Omega,\Sigma,p\right)$ be a probability measure space i.e. a measure space with  $p(\Omega)=1$ and let $k\in\N$. A $\Sigma$-measurable function $X$ from $\Omega$ to $\R^k$ is called a   (vector) \emph{random variable};  in statistical inference problems, $X$ is sometimes referred to as the given data. We write $p_X$ to denote the distribution of $X$ i.e. the measure induced by $X$ on $\left(\R^k, \mathcal{B}(\R^k)\right)$ defined by
\begin{align*}
p_X(A)=p\left(X\in A\right):=p(X^{-1}(A)),\quad \text{for all Borel set}\quad A\subseteq \R^k.
 \end{align*}
With a little abuse of terminology, $X$ is said to be an \emph{absolutely continuous random variable} if and only if $p_X\ll \mathcal L^k$  i.e. $p_X$ is absolutely continuous with respect to the Lebesgue measure $\mathcal L^k$. In this case the non-negative Radon-Nikodym derivative $f_X:=\frac{dp_X}{d \mathcal L^k}$ is called the    density  function of $X$ and it is defined through the relation
\begin{align*}
p_X(A)=\int_{A}f_X(x)\,dx,\quad \text{for all Borel set}\quad A\subseteq \R^k.
\end{align*}
$X$ is said to be a \emph{discrete random variable} if and only if there exists   a countable subset $I=(a_i)_{i\in N}$ of $\R^k$ such that $p_X\ll {\Ha^0}_{\vert I}$  i.e. $p_X$ is absolutely continuous with respect to the counting measure $\Ha^0$ on $I$. In this case the density 
 $f_X:=\frac{dp_X}{d {\Ha^0}_{\vert I}}$ is also  called the probability mass function and it is defined through the relation
\begin{align*}
p_X(A)=\sum_{i:a_i\in A}f_X(a_i),\quad \text{for all subset}\quad A\subseteq \R^k.
\end{align*}
In particular $p_X(a)=f_X(a)$, for all $a\in A$.

 Let $X:\Omega\to\R^k$ be a fixed random variable and let $n\in\N$; a random variable $Y:\Omega\to\R^n$ is called a \emph{statistic} (of the data $X$)  if it is a measurable function   of  $X$ i.e. $Y$  satisfies $Y=\phi\circ X$ for some $\left(\mathcal{B}(\R^k), \mathcal{B}(\R^n)\right)$ measurable function $\phi:\R^k\to\R^n$.

Finally, let $X_1,\dots, X_n$ be $n$ random variables where $X_i:\Omega\to\R^k$, for $i=1,\dots, n$. $X_1,\dots, X_n$ are said to be  \emph{independent} if  for every  Borel subset $A_1,\dots A_n$ of $\R^k$ and for every  $J\subseteq \{1,\dots, n\}$ one has 
\begin{align*}
p\left(\bigcap_{i\in J}X_i^{-1}(A_i)\right)=\prod_{i\in J}p_{X_i}(A_i).
\end{align*}
In this case, if every $X_i$ is absolutely continuous with density function $f_i$, then $X$ is absolutely continuous and its density function satisfies $f_X(x_1,\dots, x_n)=\prod_{i=1,\dots,n}f_i(x_i)$ for every $x_i\in\R^k$.

If, moreover, $X_1,\dots, X_n$ are identically distributed i.e. $p_{X_i}=p_{X_j}:=q$ for every $i,j$, then  $X=\left(X_1,\dots, X_n\right)$ is called a \emph{random sample} from the distribution $q$; in this case, if every $X_i$ is absolutely continuous with density $f_i=f$,  then the density function of $X$ satisfies $f_X(x_1,\dots, x_n)=\prod_{i=1,\dots,n}f(x_i)$.

\section{Area and Coarea Formulas}\label{Section AC}

In this section we provide a brief introduction to the theory of Hausdorff measures and we collect the main results about   the  Area and Coarea Formulas proved by Federer in \citep{Federer}.  For the related proofs, we refer the reader, for example,  to \citep{AFPallara, Federer-book, GiaqModica} (and references therein).

We begin with the definition of the $s$-dimensional Hausdorff measure.
Let $E\subseteq\R^n$, $\epsilon>0$ and let   $\left(B(x_i,r_i)\right)_{i\in\N}$ be  a coverings of $E$ by a countable collections of balls $B(x_i,r_i)$  whose radii  satisfy $r_i\leq \epsilon$. For each $s\geq 0$, let
\begin{align*}
\sigma_s(\epsilon)= \frac{\pi^{s/2}}{\Gamma (1 + s/2)}\inf\sum_{i\in\N} r_i^s,
\end{align*}
where the infimum is taken over all such coverings. The monotonicity of $\sigma_s$, with respect to $s$, implies that there exists the limit (finite or
infinite)
\begin{align*}
\Ha^s(E) := \lim_{\epsilon\to 0^+}\sigma_s(\epsilon).
\end{align*}
This limit is called the $s$-dimensional Hausdorff measure of E. The Hausdorff measure $\Ha^s$ satisfies Caratheodory's criterion therefore, the $\sigma$-algebra of all the  $\Ha^s$-measurable sets contains all the Borel subsets of $\R^n$ (see e.g. \citep[Proposition 2.49]{AFPallara}).

If $0<s\leq n$ is a positive
integer and $E$ is  an $s$-dimensional smooth sub-manifold of $\R^n$, then $\Ha^s(E)$ is  the $s$-dimensional volume of $E$. In particular, $\Ha^n(E) = \mathcal L^n(E)$ for Lebesgue measurable subsets $E\subseteq \R^n$ (see, e.g., \citep[Section 3.2]{Federer-book}).

When $s=0$, $\Ha^0$ coincides with the counting measure on $\R^n$ which associates to any $E\subseteq\R^n$, its number of elements $|E|$ ($+\infty$ in case the subset is infinite); in particular for every function $f:\R^n\to\R$ one has 
\begin{align*}
\int_E f(x)\,d\Ha^0(x)=\sum_{x\in E}f(x).
\end{align*}
%
%When $s=0$ $\Ha^0$ coincides with the counting measure on $\R^n$ defined by 
%\begin{align*}
%\Ha^0(E)=
%\begin{cases}
%|E|,\quad &\text{if $E$ is finite};\\
%+\infty,\quad &\text{if $E$ is infinite};
%\end{cases}
%\end{align*}

Let now  $\phi:\R^k\to\R^n$ be a locally Lipschitz map; we remark that 
%Let now  $\Omega$ be an open subset of $\R^k$ and let $\phi:\Omega\to\R^n$ be a  Lipschitz map. We remark that  it is always  possible to extend $\phi$ to the whole space preserving the Lipschitz constant  therefore, in what follows, we always suppose $\phi$ to be  globally defined; moreover
  Rademacher's Theorem assures that $\phi$ is a.e. differentiable  (see for example \citep[Proposition 2.12 and 2.14]{AFPallara}).

%%%%%%%%%%%%%%%%%%%%%%%%%%%%%%%%%%%%%%%%%%%%%%%%%%%%%%%%%%%%%%%%%%%%%%%%%

\begin{defi}[$k$-dimensional Jacobian]\label{k-Jac}
Let $k,m,n\in\N$  and let $\phi:\R^k\to\R^n$ be a locally Lipschitz map. The $m$-dimensional Jacobian of $\phi$ is defined by 
\begin{align*}
 J_m \phi(x):=\sup\left\{
 \frac{\Ha^m\Big(J\phi(x)(P)\Big)}{\Ha^m(P)}
\;:\; P\;\text{is a $m$-dimensional parallelepiped of $\R^k$}
 \right\},
 \end{align*}
 where  $J\phi(x)(P)$ is the image of $P$ under the Jacobian  matrix $J\phi(x)$ of $\phi$ which exists for a.e. $x\in\R^k$. When $\mbox{rank }(J\phi(x))\leq m$ then 
\begin{align*}
 J_m \phi(x)=\sqrt{\sum_{B}(\mbox{det\,}B)^2},
 \end{align*}
 where the sum in the last equality  runs along all $m\times m$ minors $B$ of  $J\phi(x)$ (see e.g. \citep[Section 3.6]{Morgan}).
\end{defi}
Note that $J_m\phi(x)=0$ if and only if $\mbox{rank}(J\phi(x))<m$ and that the Cauchy-Binet formula gives  in particular 
\begin{align*}
J_n \phi=\sqrt{\mbox{det}\left(J\phi\cdot J\phi^t\right)}, \qquad J_k \phi=\sqrt{\mbox{det}\left(J\phi^t\cdot J\phi\right)}.
\end{align*}

%%%%%%%%%%%%%%%%%%%%%%%%%%%%%%%%%%%%%%%%%%%%%%%%%%%%%%%%%%%%%%%%%%%%%%%%%%%%%%%%%%%%%%%
%\begin{defi}[$k$-dimensional Jacobian]\label{k-Jac}
%Let $k,m,n\in\N$ with $m\leq \min\{k,n\}$ and let $\phi:\R^k\to\R^n$ be a locally Lipschitz map. The $m$-dimensional Jacobian of $\phi$ is defined by 
%\begin{align*}
% J_m \phi=\sqrt{\sum_{B}\mbox{det}(B)^2},
% \end{align*}
% where the sum in the last equality  runs along all $m\times m$ minors $B$ of the Jacobian matrix $J\phi$ of $\phi$.
%\end{defi}
%Note that $J_m\phi=0$ if and only if $\mbox{rank}(J\phi)<m$ and that the Cauchy-Binet formula gives  in particular 
%\begin{align*}
%J_n \phi=\sqrt{\mbox{det}\left(J\phi\cdot J\phi^t\right)}, \qquad J_k \phi=\sqrt{\mbox{det}\left(J\phi^t\cdot J\phi\right)}.
%\end{align*}
In the conventional situation $k=n=m$, the above definition gives back the classic Jacobian  $J_k \phi=|\mbox{det}J\phi|$.
\bigskip

The next Theorem handles with the case $k\leq n$ and it  can be seen as a generalization of the change of variables formula \eqref{int 2} when the   invertibility and  $C^1$-regularity assumptions on $\phi$ are dropped.  For its proof we refer the reader to \citep[Theorem 3.2.11]{Federer}, \citep[Theorem 2.80]{GiaqModica}.

\begin{teo}[Area formula]\label{Area}
Let $\phi:\R^k\to\R^n$ be a locally Lipschitz map with $k\leq n$.
\begin{itemize}
\item[(i)]For any $\mathcal{L}^k$-measurable set $E\subseteq \R^k$ the multiplicity function $y\mapsto \Ha^0\left(E\cap \phi^{-1}(y)\right)$ is $\Ha^k$-measurable in $\R^n$ and
\begin{align*}
\int_E J_k\phi(x)\,dx=\int_{\R^n}\Ha^0\left(E\cap \phi^{-1}(y)\right)\,d\Ha^k(y).
\end{align*}
\item[(ii)] If $u$ is a positive measurable function,  or $u\,J_k\phi\in L^1\left(\R^k \right)$,   then
\begin{align*}
\int_{\R^k}u(x) J_k\phi(x)\,dx=\int_{\R^n}\int_{\phi^{-1}(y)}u(x)\,d\Ha^0(x)\;d\Ha^k(y)=\int_{\R^n}\sum_{x\in \phi^{-1}(y)}u(x)\;d\Ha^k(y).
\end{align*}
%\item[(iii)] If $u$ is an  $\mathcal{L}^k$-integrable function then
%\begin{align*}
%\int_{E}u(x) J_k\phi(x)\,dx=\int_{f(E)}\int_{E\cap \phi^{-1}(y)}u(x)\,d\Ha^0(x)\;d\Ha^k(y)=\int_{\phi(E)}\sum_{x\in E\cap \phi^{-1}(y)}u(x)\;d\Ha^k(y).
%\end{align*}
\end{itemize} 
\end{teo}

When $\phi: E\to \R^n$ is injective, the last formula allows the computation of the area of the  Lipschitz parametrized $k$-dimensional manifold $\phi(E)$ of $\R^n$:
\begin{align}\label{parametrized manifold}
\nonumber\Ha^k(\phi(E))&=\int_E J_k\phi\,dx=\int_E \sqrt{\mbox{det}(J\phi^t\cdot J\phi)}\,dx,\\[1ex]
\int_{\phi(E)}g(y)\,d\Ha^k(y)&=\int_E g(\phi(x))J_k\phi\,dx=\int_E g(\phi(x))\sqrt{\mbox{det}(J\phi^t\cdot J\phi)}\,dx,
\end{align}
where $g$ is any positive measurable function,  or $g\in L^1\left(\phi(E),\Ha^k \right)$.
 Note that in the particular case of a Cartesian parametrization $\phi(x) = (x,\psi(x))$ one has  $J_k\phi=\sqrt{1+\sum_B \mbox{det}(B)^2}$
where this times the sum   runs along all square minors $B$ of the Jacobian matrix $J\psi$ of $\psi$ (see \citep[page 88]{AFPallara}).

\bigskip

The next Theorem treats, conversely, the case $k> n$ and it  can be seen as a generalization of the Fubini's theorem about the reduction of  integrals. For its proof we refer the reader to \citep[Theorem 3.2.3]{Federer}, \citep[Theorem 2.86]{GiaqModica}. 
\begin{teo}[Coarea formula]\label{Coarea}
Let $\phi:\R^k\to\R^n$ be a locally Lipschitz map with $k> n$.
\begin{itemize}
\item[(i)]For any $\mathcal{L}^k$-measurable set $E\subseteq \R^k$ the  function $y\mapsto \Ha^{k-n}\left(E\cap \phi^{-1}(y)\right)$ is $\mathcal{L}^n$-measurable in $\R^n$ and
\begin{align*}
\int_E J_n\phi(x)\,dx=\int_{\R^n}\Ha^{k-n}\left(E\cap \phi^{-1}(y)\right)\,dy.
\end{align*}
\item[(ii)] If $u$ is a positive measurable function, or $uJ_n\phi\in L^1\left(\R^k\right)$, then 
\begin{align*}
\int_{\R^k}u(x) J_n\phi(x)\,dx=\int_{\R^n}\int_{\phi^{-1}(y)}u(x)\,d\Ha^{k-n}(x)\;dy.
\end{align*}
%\item[(iii)] If $u$ is an  $\mathcal{L}^k$-integrable function then
%\begin{align*}
%\int_{E}u(x) J_n\phi(x)\,dx=\int_{\phi(E)}\int_{E\cap \phi^{-1}(y)}u(x)\,d\Ha^{k-n}(x)\;dy.
%\end{align*}
\end{itemize} 
\end{teo}
When $\phi$ is an orthogonal projection 
%(e.g. $\phi(x_1,\cdots, x_k)=(x_1,\cdots, x_n)$)
(e.g. $\phi(x_1,\cdots, x_k)=(x_{i_1},\cdots, x_{i_n})$ where $\{i_1,\cdots, i_n\}\subseteq \{1,\cdots, k\}$), then $J_n\phi=1$,  the level sets of $\phi$ are $(n-k)$-planes and the last formula corresponds to Fubini's theorem.

Applying Theorem \ref{Coarea} in the particular case $n=1$, then one has  $J_1\phi(x)=|\nabla \phi(x)|$ and the formula in (ii) becomes
\begin{align}\label{1d-Coarea}
\int_{\R^k}u(x)\, |\nabla \phi(x)|\,dx=\int_{\R}\int_{\phi^{-1}(y)}u(x)\,d\Ha^{k-1}(x)\;dy.
\end{align}
In the special case $\phi(x)=|x|$, $J_n\phi(x)=1$ for every $x\neq 0$ and,  since the map sending $x\mapsto rx$ changes $\Ha^{k-1}$ by the factor $r^{k-1}$ (see e.g. \citep[Proposition 2.49]{AFPallara}),   one has  
\begin{align*}
\int_{\R^k}u(x)\,dx=\int_0^\infty\int_{|x|=r}u(x)\,d\Ha^{k-1}(x)\;dr=\int_0^\infty r^{k-1}\int_{|x|=1}u(x)\,d\Ha^{k-1}(x)\;dr.
\end{align*}

%\bigskip
%
%We end the section by stating a generalization of  the  Coarea Formula to the case where the Lebesgue measure on the right hand side of the equations in Theorem \ref{Coarea} is replaced by the Hausdorff measure $\Ha^m$, where $m=\mbox{rank }(J\phi)$. For simplicity we suppose $f$ to be $C^1$-differentiable.
%\begin{teo}\label{Coarea gen}
%Let $k,m,n\in\N$ and  let  $\phi:\R^k\to\R^n$ be a $C^1$-map. Let us fix a $\mathcal{L}^k$-measurable set $E\subseteq \R^k$  such that  $\mbox{rank }(J\phi(x))\leq m\leq k$ for every $x\in E$. The following properties hold.
%\begin{itemize}
%\item[(i)] The  function $y\mapsto \Ha^{k-m}\left(E\cap \phi^{-1}(y)\right)$ is $\Ha^m$-measurable in $\R^n$ and
%\begin{align*}
%\int_E J_m\phi(x)\,dx=\int_{\R^n}\Ha^{k-m}\left(E\cap \phi^{-1}(y)\right)\,d\Ha^m(y).
%\end{align*}
%\item[(ii)] If $u$ is a positive measurable function, or $uJ_m\phi\in L^1\left(E\right)$, then 
%\begin{align*}
%\int_{E}u(x) J_m\phi(x)\,dx=\int_{\R^n}\int_{\phi^{-1}(y)\cap E}u(x)\,d\Ha^{k-m}(x)\;d\Ha^m( y).
%\end{align*}
%\end{itemize} 
%\end{teo}
% {\sc{Proof.}} The proof is a consequence of  \citep[Theorem 5.1, Theorem 5.2]{KristensenABridge}.\qed

\bigskip

We end the section by stating a generalization of  the  Coarea Formula to the case where the Lebesgue measure on the right hand side of the equations in Theorem \ref{Coarea} is replaced by the Hausdorff measure $\Ha^m$, where $m:=\max\{\mbox{rank }(J\phi(x)):x\in\R^k\}$. For simplicity we suppose $f$ to be $C^1$-differentiable.
\begin{teo}\label{Coarea gen}
	Let $k,n\in\N$ and  let  $\phi:\R^k\to\R^n$ be a $C^1$-map and let   $m:=\max\{\mbox{rank }(J\phi(x)):x\in\R^k\}$ . The following properties hold.
	\begin{itemize}
		\item[(i)] For every  $\mathcal{L}^k$-measurable set $E\subseteq \R^k$ the  function $y\mapsto \Ha^{k-m}\left(E\cap \phi^{-1}(y)\right)$ is $\Ha^m$-measurable in $\R^n$ and one has
		\begin{align*}
			\int_E J_m\phi(x)\,dx=\int_{\R^n}\Ha^{k-m}\left(E\cap \phi^{-1}(y)\right)\,d\Ha^m(y).
		\end{align*}
		\item[(ii)] If $u$ is a positive measurable function, or $uJ_m\phi\in L^1\left(E\right)$, then 
		\begin{align*}
			\int_{E}u(x) J_m\phi(x)\,dx=\int_{\R^n}\int_{\phi^{-1}(y)\cap E}u(x)\,d\Ha^{k-m}(x)\;d\Ha^m( y).
		\end{align*}
	\end{itemize} 
\end{teo}
{\sc{Proof.}} The proof is a consequence of  \citep[Theorem 5.1, Theorem 5.2]{KristensenABridge}.\qed

\bigskip
The next Remark clarifies some  positivity properties about   $k$-dimensional Jacobians. It  can be seen as a generalization of Sard's Theorem (see also \citep[Lemma 2.73, Lemma 2.96, Remark 2.97]{AFPallara} and \citep[Theorem 1.1]{KristensenABridge}).

\begin{os}\label{oss J>0}
Let  $\phi:\R^k\to \R^n$ be  a locally Lipschitz map and let us suppose, preliminarily, that  $J\phi(x)$ has maximum rank for a.e. $x\in\R^k$.
\begin{itemize}
\item[(i)] if $k\leq n$, then using (i) of  Theorem \ref{Area}  with $E:=\{x\in\R^k\,:\,J_k\phi(x)=0\}$ we get
\begin{align*}
\int_{\R^n}\Ha^0\left(E\cap \phi^{-1}(y)\right)\,d\Ha^k(y)=0.
\end{align*}
This implies $\Ha^0\left(E\cap \phi^{-1}(y)\right)=0$ (i.e. $\phi(E)\cap \{y\}=\emptyset $)  for  $\Ha^k$-a.e. $y\in\R^n$. This yields  $\Ha^k\left(\phi(E)\right)$=0 and it  implies, in particular,  that $J_k\phi>0$ on $\phi^{-1}(y)$ for  $\Ha^k$-a.e. $y\in\R^n$.
\item[(ii)] if $k\geq n$, then using (i) of  Theorem \ref{Coarea}  with $E:=\{x\in\R^k\,:\,J_n\phi(x)=0\}$ we get
\begin{align*}
\int_{\R^n}\Ha^{k-n}\left(E\cap \phi^{-1}(y)\right)\,dy=0.
\end{align*}
This yields $\Ha^{k-n}\left(E\cap \phi^{-1}(y)\right)=0$   for  a.e. $y\in\R^n$ and it implies, in particular,  that $J_n\phi>0$ $\Ha^{k-n}$-a.e on $\phi^{-1}(y)$ for   a.e. $y\in\R^n$.
\end{itemize}
\medskip
Let us suppose, now, $\phi:\R^k\to\R^n$ to be  $C^1$ and let   $m:=\max\{\mbox{rank }(J\phi(x)):x\in\R^k\}$.   Then  setting  $E:=\{x\in\R^k\,:\,J_m\phi(x)=0\}$ and  using Theorem \ref{Coarea gen} we get
\begin{align*}
\int_{\R^n}\Ha^{k-m}\left(E\cap \phi^{-1}(y)\right)\,d\Ha^m(y)=0.
\end{align*}
This implies that $J_m\phi>0$ $\Ha^{k-m}$-a.e on $\phi^{-1}(y)$ for   $\Ha^m$-a.e. $y\in\R^n$.

%Moreover, using Theorem \ref{Coarea gen}, if $\phi:\R^k\to\R^n$ is  $C^1$ and   $\mbox{rank }(J\phi(x))\leq m\leq k$ for every $x\in\R^k$, then  setting  $E:=\{x\in\R^k\,:\,J_m\phi(x)=0\}$ we get
%\begin{align*}
%\int_{\R^n}\Ha^{k-m}\left(E\cap \phi^{-1}(y)\right)\,d\Ha^m(y)=0.
%\end{align*}
%This implies that $J_m\phi>0$ $\Ha^{k-m}$-a.e on $\phi^{-1}(y)$ for   $\Ha^m$-a.e. $y\in\phi(\R^k)\subseteq\R^n$.
\end{os}
\section{Sample Distribution Theory}\label{Section Distribution}
Let $\left(\Omega,\Sigma,p\right)$ be a probability measure space, let $k\in\N$ and let $X:\Omega\to \R^k$ be an absolutely  continuous random variable. Let $Y:=\phi\circ X$ be a statistic, where $\phi:\R^k\to \R^n$ is a measurable map and  $k\in\N$. In this section we prove that when  $\phi$ is  locally Lipschitz then  
%$Y$ is a continuous random variable  and we  compute its probability density function in terms of an integral involving the  density function $f_X$ of the data $X$.
the probability measure $p_Y$ induced by $Y$ has a density function, with respect to some Hausdorff measure $\Ha^m$ on $\phi(\R^k)\subseteq\R^n$, which can be computed explicitly in terms of an integral involving the  density function $f_X$ of  $X$.
 We recall preliminarily that the  Radon-Nikodym derivative of a measure  is uniquely defined up to   zero measure   sets: since, by  definition,  $p_Y$ is concentrated on $\phi(\R^k)$, in what follows we can always set  $f_Y(y)=0$ for  $y\notin\phi(\R^k)$.

We start with the case $k\leq n$.

%\begin{os}
%Let $\mu,\nu$ two measures defined on a measure space $\left(\Omega,\Sigma\right)$ such that $\nu\ll\mu$ and let $\frac{d\nu}{d\mu}$ be the Radon-Nikodym derivative of $\nu$ with respect to $\mu$.  $\frac{d\nu}{d\mu}$ is uniquely defined up to   zero $\mu$-measure   sets.  In particular if $\mu=\mathcal{L}^n$ we can suppose $\frac{d\nu}{d\mu}=0$ on every hyperplanes of $\R^n$ of dimension $k<n$. In particular it is sufficient to define  $\frac{d\nu}{d\mu}$ on the set  $\{(x_1,\dots,x_n)\in\R^n\;:\;x_i\neq x_j,\, i\neq j\}$.
%\end{os}

\begin{teo}\label{teo k<n}
Let $\left(\Omega,\Sigma,p\right)$ be a probability measure space, let $k,n \in\N$ with $k\leq n$ and let $X:\Omega\to \R^k$ be an absolutely  continuous random variable with probability density function $f_X$. If $\phi:\R^k\to \R^n$ is a locally Lipschitz map  such that $\mbox{rank }(J\phi)=k$ a.e.,  then the probability measure induced by the  statistic $Y:=\phi\circ X$ is absolutely continuous with respect to the Hausdorff measure $\Ha^k$ on $\R^n$ i.e.  $p_{Y}\ll \Ha^k$. Its Radon-Nykodym derivative $\frac{dp_y}{d\mathcal \Ha^k}$ is defined through the relation 
\begin{align*}
p_Y(A)=p\left(Y^{-1}(A)\right)=\int_{A}\frac{dp_Y}{d\mathcal \Ha^k}(y)\,d\Ha^{k}(y),\quad \text{for all Borel  subset}\quad A\subseteq \R^n.
\end{align*}
%\begin{align*}
%p_Y(A)=p\left(Y^{-1}(A)\right)=\int_{A}\frac{dp_Y}{d\mathcal \Ha^k}(y)\,d\Ha^{k}(y),\quad \text{for all $\Ha^k$-measurable set}\quad A\subseteq \R^n
%\end{align*}
It satisfies
%\begin{align*}
%\frac{dp_Y}{d\mathcal \Ha^k}(y)=\begin{cases}
%\int_{\phi^{-1}(y)}f_X(x)\frac{1}{J_k\phi(x)}\,d\Ha^0(x)=\sum_{\phi(x)=y}f_X(x)\frac{1}{J_k\phi(x)}, &\quad \text{for $\Ha^k$-a.e.}\quad y\in\phi(R^k);\\[3ex]
%0, &\quad \text{otherwise}.
%\end{cases}
%\end{align*}
\begin{align*}
\frac{dp_Y}{d\mathcal \Ha^k}(y)&=
\int_{\phi^{-1}(y)}f_X(x)\frac{1}{J_k\phi(x)}\,d\Ha^0(x)=\sum_{\phi(x)=y}f_X(x)\frac{1}{J_k\phi(x)}, &\quad \text{for $\Ha^k$-a.e.}\quad y\in\phi(R^k)
\end{align*}
and it is $0$ otherwise.

% Moreover let us suppose that there exists a countable disjoint covering $\bigcup_{i\in\N} E_i$ of $\Ha^k$-a.e point  of $\R^n$ (i.e. the set of point of $\R^n$ which is not covered has $\Ha^k$-measure zero)  such that on each measurable subset   $E_i$,  $\phi_i:=\phi_{\vert E_i}$ is $\Ha^k$-a.e injective. Then 
Moreover let us suppose that there exists a countable disjoint covering $\bigcup_{i\in\N} E_i$ of a.e point  of $\R^k$ (i.e. the set of points of $\R^k$ which are not covered has $\mathcal{L}^k$-measure zero)  such that on each measurable subset   $E_i$, the restriction map  $\phi_i:=\phi_{\vert E_i}$ is a.e injective. Then 
\begin{align*}
\frac{p_Y}{d\mathcal \Ha^k}(y)=\sum_{i\in \N}f_X\left(\phi_i^{-1}(y)\right)\frac{1}{J_k\phi\left(\phi_i^{-1}(y)\right)},\quad\text{for $\Ha^k$-a.e.}\quad y\in\phi(R^k).
\end{align*}
\end{teo}
{\sc{Proof.}} Let $A\subseteq \R^n$ be a Borel set. Recalling Definition \ref{k-Jac} and Remark \ref{oss J>0},   $J_k\phi>0$ on $\phi^{-1}(y)$ for  $\Ha^k$-a.e. $y\in\R^n$. Then using the Area Formula of Theorem \ref{Area}  one has
\begin{align*}
p_Y(A)&=p\left(Y^{-1}(A)\right)=p\left(X^{-1}\left(\phi^{-1}(A)\right)\right)=\int_{\phi^{-1}(A)}f_X(x)\,dx\\[1ex]
&=\int_{\R^n}\int_{\phi^{-1}(y)\cap \phi^{-1}(A)}f_X(x)\frac{1}{J_k\phi(x)}\,d\Ha^0(x)\;d\Ha^k(y)\\[1ex]
&=
\int_{A}\int_{\phi^{-1}(y)}f_X(x)\frac{1}{J_k\phi(x)}\,d\Ha^0(x)\;d\Ha^k(y)=\int_{A}\sum_{x\in \phi^{-1}(y)}f_X(x)\frac{1}{J_k\phi(x)}\;d\Ha^k(y).
\end{align*}
This proved the first  required claim. The second assertion follows after observing that, under the given hypothesis, $\phi^{-1}(y)=\bigcup_i \{\phi_i^{-1}(y)\}$ for every $y\in \phi(\R^k)$.
\\\qed
\begin{os}
We note that if $k<n$, then $H^k$ is not $\sigma$-finite on $\R^n$, so we cannot directly use  the Radon-Nikódym theorem in order to deduce from $p_{Y}\ll \Ha^k$ the existence of $\frac{dp_Y}{d\mathcal \Ha^k}(y)$. Nevertheless in this case $p_Y$ is concentrated on $\phi(\R^k)$ which has $\sigma$-finite $\Ha^k$-measure and therefore the Radon-Nikódym theorem applies using  $p_{Y}\ll \Ha^k_{\vert \phi(\R^k)}$. 

Indeed if $\R^k\subseteq \bigcup_{i\in\N} E_i$, where each $E_i$ is a Borel subset of $\R^k$ such that 
$\mathcal{L}^k\left(E_i\right)<\infty$, then 
$\phi\left(\R^k\right)\subseteq \bigcup_{i\in\N} \phi\left(E_i\right)$ 
and from \citep[Proposition 2.49]{AFPallara} one has  $\Ha^k\left(\phi\left(E_i\right)\right)<\mbox{Lip}(\phi)^k \mathcal{L}^k\left(E_i\right)<\infty$.
\end{os}
When  $k=n$, then recalling that $\Ha^k=\mathcal{L}^k$, the previous theorem implies that   $p_Y \ll  \mathcal{L}^k$.
\begin{cor}\label{teo k=n}
Let $\left(\Omega,\Sigma,p\right)$ be a probability measure space, let $k\in\N$ and let $X:\Omega\to \R^k$ be an absolutely  continuous random variable with probability density function $f_X$. If $\phi:\R^k\to \R^k$ is a locally  Lipschitz map such that $\mbox{rank }(J\phi)=k$ a.e.,  then the statistic $Y:=\phi\circ X$ is an absolutely  continuous random variable  and its probability density function  $f_Y$  satisfies
\begin{align*}
f_Y(y)=\sum_{\phi(x)=y}f_X(x)\frac{1}{J_k\phi(x)}, &\quad \text{for a.e.}\quad y\in\phi(R^k)
\end{align*}
and it is $0$ otherwise.  Moreover let us suppose that there exists a countable disjoint covering $\bigcup_{i\in\N} E_i$ of a.e point  of $\R^k$  such that on each measurable subset   $E_i$,  $\phi_i:=\phi_{\vert E_i}$ is a.e injective. Then 
\begin{align*}
f_Y(y)=\sum_{i\in \N}f_X\left(\phi_i^{-1}(y)\right)\frac{1}{J_k\phi\left(\phi_i^{-1}(y)\right)},\quad \text{for a.e.}\quad y\in\phi(R^k).
\end{align*}
\end{cor}
\bigskip

Let us now consider  the case $k>n$.

\begin{teo}\label{teo k>n}
Let $\left(\Omega,\Sigma,p\right)$ be a probability measure space, let $k,n \in\N$ with $k\geq n$ and let $X:\Omega\to \R^k$ be an absolutely  continuous random variable with probability density function $f_X$. If $\phi:\R^k\to \R^n$ is a locally Lipschitz map such that $\mbox{rank }(J\phi)=n$ a.e.,  then the statistic $Y:=\phi\circ X$ is an absolutely  continuous random variable (i.e.  $p_{Y}\ll \mathcal{L}^n$)  and its probability density function  $f_Y$  satisfies  
%\begin{align*}
%f_Y(y)=
%\begin{cases}
%\int_{\phi^{-1}(y)}f_X(x)\frac{1}{J_nf(x)}\,d\Ha^{k-n}(x), &\quad \text{for a.e.}\quad y\in\phi(R^k);\\[3ex]
%0, &\quad \text{otherwise}.
%\end{cases}
%\end{align*} 
\begin{align*}
f_Y(y)=
\int_{\phi^{-1}(y)}f_X(x)\frac{1}{J_n\phi(x)}\,d\Ha^{k-n}(x), &\quad \text{for a.e.}\quad y\in\phi(R^k)
\end{align*} 
and it is $0$ otherwise.
\end{teo}
{\sc{Proof.}} The case $k=n$ is the result of  the previous Corollary. Let us suppose $k>n$. Recalling Definition \ref{k-Jac} and Remark \ref{oss J>0},  $J_n\phi>0$ $\Ha^{k-n}$-a.e on $\phi^{-1}(y)$ for   a.e. $y\in\R^n$. Let $A\subseteq \R^n$ be a Borel set. Then using the Coarea Formula of Theorem \ref{Coarea}  one has
\begin{align*}
p_Y(A)&=p\left(Y^{-1}(A)\right)=p\left(X^{-1}\left(\phi^{-1}(A)\right)\right)\\[1ex]
&=\int_{\phi^{-1}(A)}f_X(x)\,dx=\int_{\R^n}\int_{\phi^{-1}(y)\cap \phi^{-1}(A)}f_X(x)\frac{1}{J_n\phi(x)}\,d\Ha^{k-n}(x)\;dy\\[1ex]
&=
\int_{A}\int_{\phi^{-1}(y)}f_X(x)\frac{1}{J_n\phi(x)}\,d\Ha^{k-n}(x)\;dy.
\end{align*}
This proved the required claim.\\\qed

For  the reader's convenience we enlighten in the following corollary the particular case $n=1$ which is very useful in the applications and which follows from formula \eqref{1d-Coarea}.   

\begin{cor}\label{teo k>1}
Let $\left(\Omega,\Sigma,p\right)$ be a probability measure space, let $k\in\N$  and let $X:\Omega\to \R^k$ be an absolutely  continuous random variable with probability density function $f_X$. If $\phi:\R^k\to \R$ is a locally Lipschitz map such that $|\nabla \phi|>0$ a.e.,  then the statistic $Y:=\phi\circ X$ is an absolutely  continuous random variable   and has probability density function 
%\begin{align*}
%f_Y(y)=
%\begin{cases}
%\int_{\phi^{-1}(y)}f_X(x)\frac{1}{J_nf(x)}\,d\Ha^{k-n}(x), &\quad \text{for a.e.}\quad y\in\phi(R^k);\\[3ex]
%0, &\quad \text{otherwise}.
%\end{cases}
%\end{align*} 
\begin{align*}
f_Y(y)=
\int_{\phi^{-1}(y)}f_X(x)\frac{1}{|\nabla\phi(x)|}\,d\Ha^{k-1}(x), &\quad \text{for a.e.}\quad y\in\phi(R^k).
\end{align*} 
\end{cor}
\medskip
\begin{os}The assumptions $p(\Omega)=1$ was never used in the proof of the Theorems \ref{teo k<n} and \ref{teo k>n}. Indeed  analogous results hold  with $p_X$ replaced by an absolutely continuous measure on $\R^k$. More precisely let $\mu$ be  a measure defined on $\left(\R^k,\mathcal{B}\right)$  such that $\mu \ll \mathcal L^k$ and let $\frac{d\mu}{d\mathcal L^k}$ its Radon-Nykodym derivative. Let  $\phi:\R^k\to\R^n$ be a locally Lipschitz map whose Jacobian matrix $J\phi$ has a.e. maximum rank.
\begin{itemize}
\item[(i)] If $k\leq n$ then $\mu \phi^{-1} \ll  \Ha^k$ and for $\Ha^k$-a.e. $y\in \phi(\R^k)$one has
\begin{align*}
\frac{d\mu \phi^{-1}}{d\mathcal \Ha^k}(y)=\int_{\phi^{-1}(y)}\frac{d\mu}{d\mathcal L^k}(x)\frac{1}{J_k\phi(x)}\,d\Ha^0(x)=\sum_{\phi(x)=y}\frac{d\mu}{d\mathcal L^k}(x)\frac{1}{J_k\phi(x)}.
\end{align*}
%(In particular, if $k=n$, then $\mu f^{-1} \ll  \mathcal{L}^n$). Moreover if there exists a countable  covering $\R^n=\bigcup_{i\in\N} E_i$ such that on each measurable subset   $E_i$,  $f_i:=f_{\vert E_i}$ is a.e injective then 
%\begin{align*}
%\frac{d\mu f^{-1}}{d\mathcal \Ha^k}(y)=\sum_{i\in \N}\frac{d\mu}{d\mathcal L^k}\left(f_i^{-1}(y)\right)\frac{1}{J_kf\left(f_i^{-1}(y)\right)}.
%\end{align*}
\item[(ii)] If $k> n$ then $\mu \phi^{-1} \ll  \mathcal{L}^n$ and for a.e. $y\in\R^n$ one has
\begin{align*}
\frac{d\mu \phi^{-1}}{d\mathcal{L}^n}(y)=\int_{\phi^{-1}(y)}\frac{d\mu}{d\mathcal{L}^k}(x)\frac{1}{J_n\phi(x)}\,d\Ha^{k-n}(x).
\end{align*}
\end{itemize} 
\end{os}
\bigskip

We end the section by applying Theorem \ref{Coarea gen}  in order to  extend Theorem \ref{teo k>n}  to the case of a  $C^1$-map $\phi$ whose Jacobian could possibly have not maximum  rank.    In this case, setting $m:=\max\{\mbox{rank }(J\phi(x)):x\in\R^k\}$, the induced probability $p_Y$  has a density function $f_Y$ with respect to the  Hausdorff measure $\Ha^m$ on $\phi(\R^k)\subseteq\R^n$.

\begin{teo}\label{teo Hausdorff}
Let $\left(\Omega,\Sigma,p\right)$ be a probability measure space, let $k,n \in\N$ 
%with $k\geq n$
 and let $X:\Omega\to \R^k$ be an absolutely  continuous random variable with probability density function $f_X$. Let $\phi:\R^k\to \R^n$ be a $C^1$-map  and let $m:=\max\{\mbox{rank }(J\phi(x)):x\in\R^k\}$. Then the induced probability measure $p_Y$ of the statistic $Y:=\phi\circ X$ has a density function $f_Y$ with respect to the  Hausdorff measure $\Ha^m$ 
 % on $\phi(\R^k)\subseteq\R^n$ 
 which  satisfies  
%\begin{align*}
%f_Y(y)=
%\begin{cases}
%\int_{\phi^{-1}(y)}f_X(x)\frac{1}{J_nf(x)}\,d\Ha^{k-n}(x), &\quad \text{for a.e.}\quad y\in\phi(R^k);\\[3ex]
%0, &\quad \text{otherwise}.
%\end{cases}
%\end{align*} 
\begin{align*}
f_Y(y)=
\int_{\phi^{-1}(y)}f_X(x)\frac{1}{J_m\phi(x)}\,d\Ha^{k-m}(x), &\quad \text{for $\Ha^m$-a.e.}\quad y\in\phi(R^k)
\end{align*} 
and it is $0$ otherwise.
\end{teo}
%\begin{teo}\label{teo Hausdorff}
%Let $\left(\Omega,\Sigma,p\right)$ be a probability measure space, let $k,n \in\N$ 
%%with $k\geq n$
% and let $X:\Omega\to \R^k$ be a continuous random variable with probability density function $f_X$. Let $\phi:\R^k\to \R^n$ be a $C^1$-map such that $\mbox{rank }(J\phi(x))= m$ for every $x\in\R^k$. Then the induced probability measure $p_Y$ of the statistic $Y:=\phi\circ X$ has a density function $f_Y$ with respect to the  Hausdorff measure $\Ha^m$
% % on $\phi(\R^k)\subseteq\R^n$ 
% which  satisfies  
%%\begin{align*}
%%f_Y(y)=
%%\begin{cases}
%%\int_{\phi^{-1}(y)}f_X(x)\frac{1}{J_nf(x)}\,d\Ha^{k-n}(x), &\quad \text{for a.e.}\quad y\in\phi(R^k);\\[3ex]
%%0, &\quad \text{otherwise}.
%%\end{cases}
%%\end{align*} 
%\begin{align*}
%f_Y(y)=
%\int_{\phi^{-1}(y)}f_X(x)\frac{1}{J_m\phi(x)}\,d\Ha^{k-m}(x), &\quad \text{for $\Ha^m$-a.e.}\quad y\in\phi(R^k)
%\end{align*} 
%and it is $0$ otherwise.
%\end{teo}
{\sc{Proof.}} Recalling Definition \ref{k-Jac} and Remark \ref{oss J>0},  $J_m\phi>0$ $\Ha^{k-m}$-a.e on $\phi^{-1}(y)$ for   $\Ha^m$-a.e. $y\in\R^n$. Let $A\subseteq \R^n$ be a Borel set. Then using  Theorem \ref{Coarea gen}  one has
\begin{align*}
p_Y(A)&=p\left(Y^{-1}(A)\right)=p\left(X^{-1}\left(\phi^{-1}(A)\right)\right)\\[1ex]
&=\int_{\phi^{-1}(A)}f_X(x)\,dx=\int_{\R^n}\int_{\phi^{-1}(y)\cap \phi^{-1}(A)}f_X(x)\frac{1}{J_m\phi(x)}\,d\Ha^{k-m}(x)\;d\Ha^m(y)\\[1ex]
&=
\int_{A}\int_{\phi^{-1}(y)}f_X(x)\frac{1}{J_m\phi(x)}\,d\Ha^{k-m}(x)\;d\Ha^m(y).
\end{align*}
This proves the required claim.\\\qed

%\textcolor{red}{\begin{itemize}
%\item[(i)]$p_Y$ is concentrated on $\phi(\R^k)\subseteq \R^n$ which by Sard Theorem (Wikipedia) is $m$-dimensional?
%\item[(ii)] $\Ha^m_{\vert \phi(\R^k)}$ is $\sigma$-finite? (in positive case we can delete definition of measures with density in Preliminaries. 
%\end{itemize}
%\begin{prop}
%Let $\left(\Omega,\Sigma,p\right)$ be a probability measure space, let $k\in\N$ and let $X:\Omega\to \R^k$ be a continuous random variable. Let $Y:=\phi\circ X$ be a statistic, with $\phi:\R^k\to \R^n$, $k\in\N$. The following properties hold:
%\begin{itemize}
%\item[(i)] The probability measure $p_Y$ induced by $Y$  is concentrated on the Borel subset $\phi(\R^k)\subseteq\R^n$ i.e. $P_Y\left(\phi(\R^k)^c\right)=0$.
%\item[(ii)] $\Ha^m_{\vert \phi(\R^k)}$ is $\sigma$-finite
%\end{itemize}  
%\end{prop}
%{\sc{Proof.}} \citep[Proposition 2.49]{AFPallara}
%}
\section{Further generalizations}\label{Section Generalization}

In this section we  briefly expose, for the interested reader,  a further  generalization of  Theorem \ref{teo k>n} which covers  the case of  probabilities $p_X$  concentrated on some $m$-dimensional subset $E\subseteq\R^k$. This includes, for example, the cases of a random variable $X$ which is uniformly distributed over a generic subset $E$ and  of random variables taking values on $m$-dimensional manifolds.  Statistical analysis on   manifolds 
	has many applications in directional and axial statistics, morphometrics, medical
	diagnostics, machine vision and image analysis (see e.g \citep{Bhattacharya} and references there in).  Amongst the many important applications, those arising, for example, from the analysis of  data on torus  play also a fundamental role in molecular biology in the study of the Protein Folding Problem. 

Although the following Theorem is valid for   countably $m$-rectifiable sets (see \citep[Definition 5.4.1, Lemma 5.4.2]{Krantz-Parks}), we suppose for simplicity $E$ to be  an $m$-dimensional sub-manifold of $\R^k$. 

We state first the Area and Coarea formula relative to   sub-manifolds of $\R^k$. If $J\phi^E$ is the  tangential Jacobian matrix of $\phi$ with respect to $E$ (see \citep[Formula 11.1]{Maggi}), the $k$-tangential Jacobian $J_k^E\phi$ is defined as in Definition \ref{k-Jac}. For a rigorous introduction to tangential Jacobians  as well as for all the other details we refer the reader to \citep[Chapter 3]{Federer}, \citep[Section 5.3]{Krantz-Parks},  \citep[Section 11.1]{Maggi}.

\begin{teo}
Let $\phi:\R^k\to\R^n$ be a $C^1$-map  and let $E\subseteq\R^k$ an $m$-dimensional manifold. The following properties hold.
\begin{itemize}
\item[(i)]If $m\leq n$ and $u$ is a positive measurable function, or $uJ_m^E\phi\in L^1\left(E,\Ha^m\right)$, one has 
\begin{align*}
\int_{E}u\, J_m^E\phi\,d\Ha^m(x)=\int_{\R^n}\int_{E\cap \phi^{-1}(y)}u\,d\Ha^{0}(x)\,d\Ha^{m}(y).
\end{align*}
\item[(ii)]
If $m\geq n$ and  $u$ is a positive measurable function, or $uJ_n^E\phi\in L^1\left(E,\Ha^m\right)$,one has
\begin{align*}
\int_{E}u\, J_n^E\phi\,d\Ha^m(x)=\int_{\R^n}\int_{E\cap \phi^{-1}(y)}u\,d\Ha^{m-n}(x)\,dy.
\end{align*}
\end{itemize}
\end{teo}
{\sc{Proof.}} See \citep[Theorem 2.91, 2.93]{Federer} and \citep[Theorem 2.91, Theorem 2.93]{AFPallara}.\qed
\bigskip
The same methods of proof used in Section \ref{Section Distribution},  yield finally  the next result. Note that every  $m$-dimensional manifold $E\subseteq\R^k$ has $\Ha^m$-$\sigma$-finite measure.
\begin{teo}
Let $\left(\Omega,\Sigma,p\right)$ be a probability measure space, let $k,n \in\N$ and let $E\subseteq\R^k$ an $m$-dimensional manifold. Let $X:\Omega\to E$ be a random variable having a  probability density function $f_X$ with respect to the  Hausdorff measure $\Ha^m_{\vert E}$ on $E$. Let  $\phi:\R^k\to \R^n$ be  a $C^1$-map whose tangential Jacobian  $J^E\phi(x)$ has maximum rank at any point $x\in E$.  The following properties hold.
\begin{itemize}
\item[(i)] If $m\leq n$ then the probability measure $p_Y$ induced by the  statistic $Y:=\phi\circ X$ is absolutely continuous    with respect to the Hausdorff measure $\Ha^m_{\vert\phi(E)}$ on $\phi(E)$  and its density function $f_Y$  satisfies
\begin{align*}
f_Y(y)&=
\int_{\phi^{-1}(y)}f_X(x)\frac{1}{J_m^E\phi(x)}\,d\Ha^0(x)=\sum_{\phi(x)=y}f_X(x)\frac{1}{J_m^E\phi(x)}, &\quad \text{for $\Ha^m$-a.e.}\quad y\in\phi(E).
\end{align*}
\item[(ii)] If $m\geq n$ then  the statistic $Y:=\phi\circ X$ is an absolutely  continuous random variable (i.e.  $p_{Y}\ll \mathcal{L}^n$)  and its probability density function  $f_Y$  satisfies  
\begin{align*}
f_Y(y)=
\int_{\phi^{-1}(y)}f_X(x)\frac{1}{J_n^E\phi(x)}\,d\Ha^{m-n}(x), &\quad \text{for a.e.}\quad y\in\phi(E).
\end{align*} 
\end{itemize}
%
%
%\begin{itemize}
%\item[(i)] If $m\leq n$ then the probability measure induced by the  statistic $Y:=\phi\circ X$  has a density function $f_Y$ with respect to the Hausdorff measure $\Ha^m_{\vert\phi(E)}$ on $\phi(E)$ 
%\begin{align*}
%f_Y(y)&=
%\int_{\phi^{-1}(y)}f_X(x)\frac{1}{J_m^E\phi(x)}\,d\Ha^0(x)=\sum_{\phi(x)=y}f_X(x)\frac{1}{J_m^E\phi(x)}, &\quad \text{for $\Ha^m$-a.e.}\quad y\in\phi(E).
%\end{align*}
%\item[(ii)] If $m\geq n$ then  the statistic $Y:=\phi\circ X$ is a continuous random variable (i.e.  $p_{Y}\ll \mathcal{L}^n$)  and its probability density function  $f_Y$  satisfies  
%\begin{align*}
%f_Y(y)=
%\int_{\phi^{-1}(y)}f_X(x)\frac{1}{J_n^E\phi(x)}\,d\Ha^{m-n}(x), &\quad \text{for a.e.}\quad y\in\phi(E).
%\end{align*} 
%\end{itemize}
\end{teo}

\section{Some applications}\label{Section apll}

In this section we apply the results of the previous sections in order  to compute the density functions of some distributions in some   cases of relevant interest.
\subsection{First examples}
In these first examples we provide a density formula for random variables which are algebraic manipulations of absolutely continuous random variables. The first example, in particular, is used in Proposition \ref{section chi squared} in order to find the density of the chi-squared distribution.
\begin{esem}[Square function]\label{Square function}
Let $\left(\Omega,\Sigma,p\right)$ be a probability measure space and  let $X:\Omega\to \R$ be an absolutely  continuous random variable with probability density function $f_X$. We employ  Corollary \ref{teo k=n} with  \begin{align*}
\phi:\R\to [0,\infty[,\quad \phi(t)=t^2,\quad J_1(t)=2|t|.
\end{align*}
Then  the statistic $Y:=X^2$ is an absolutely  continuous random variable  and its probability density function  $f_Y$  satisfies  
\begin{align*}
f_Y(y)=\frac{f_X(\sqrt y)+f_X\left(-\sqrt y\right)}{2\sqrt y},\quad \text{for any}\quad y>0.
\end{align*}
More generally let  $k\in\N$ and let  $X=(X_1,\dots, X_k):\Omega\to \R^k$ be an absolutely  continuous (vector valued) random variable with probability density function $f_X$. Then employing   Theorem \ref{teo k>n} with 
\begin{align*}
\phi:\R^k\to [0,\infty[,\quad \phi(x)=\|x\|^2=x_1^2+\cdots+x_k^2,\quad J_1(x)=2\|x\|,
\end{align*}
 we get that the statistic $Y:=\|X\|^2=X_1^2+\cdots+X_k^2$ is an absolutely  continuous random variable whose probability density function satisfies  
\begin{align*}
f_Y(y)&=\int_{\|x\|^2=y}\frac{f_X(x)}{2\|x\|}\,d\Ha^{k-1}(x)=\frac{1}{2\sqrt y}\int_{\|x\|=\sqrt y}f_X(x)\,d\Ha^{k-1}(x),\quad \text{for any}\quad y>0.
\end{align*}
\end{esem}

%\begin{esem}[Affine transformations]
%Let $\left(\Omega,\Sigma,p\right)$ be a probability measure space, let $k\in\N$ and  let $X:\Omega\to \R$ be a continuous random variable with probability density function $f_X$. We employ  Corollary \ref{teo k=n} with  \begin{align*}
%\phi:\R^k\to \R^k,\quad \phi(t)=At+x_0,\quad J_k(t)=|\mbox{det } A|,
%\end{align*}
%where $A\in\R^{k\times k}$ is a not-singular matrix and $x_0\in\R^k$. 
%Then  the statistic $Y:=AX+x_0$ is a continuous random variable  and its probability density function  $f_Y$  satisfies  
%
%\begin{align*}
%f_Y(y)=
%\int_{Ax+x_0=y}\frac{f_X(x)}{|\mbox{det } A|}\,d\Ha^{0}(x)=\frac{f_X(A^{-1}(y-x_0))}{|\mbox{det } A|}, &\quad \text{for a.e.}\quad y\in \R^k.
%\end{align*} 
%\end{esem}

\begin{esem}[Affine transformations]\label{Affine transformations}
Let $\left(\Omega,\Sigma,p\right)$ be a probability measure space, let $k\in\N$ and  let $X:\Omega\to \R^k$ be an absolutely  continuous random variable with probability density function $f_X$.  Let us consider the affine transformation 
\begin{align*}
\phi:\R^k\to \R^n,\quad \phi(x)=Ax+y_0,
\end{align*}
where $A\in\R^{n\times k}$,  $\mbox{rank}(A)=m$ and $y_0\in\R^n$. Recalling Definition \ref{k-Jac}, the $m$-dimensional Jacobian of $\phi$ is given by
\begin{align*}
	J_m\phi(x)=\sqrt{\sum_{B}(\mbox{det\,}B)^2}=:A_m
\end{align*}
where the sum   runs along all $m\times m$ minors $B$ of  $A$.
 Then, using Theorem \ref{teo Hausdorff},  the induced probability measure $p_Y$ of the statistic $Y:=AX+y_0$ has a density function $f_Y$ with respect to the  Hausdorff measure $\Ha^m$ on the $m$-dimensional hyper-surface $\phi\left(\R^k\right)=\{y=Ax+y_0:x\in\R^k\}$ which  satisfies   
\begin{align*}
f_Y(y)&=
\frac{1}{A_m}\int_{Ax+y_0=y}f_X(x)\,d\Ha^{k-m}(x)
\\[1ex]&=\frac{1}{A_m}\int_{\mbox{Ker}(A)+x_y}f_X(x)\,d\Ha^{k-m}(x)
, \quad \text{for}\quad y\in\phi(R^k).
\end{align*} 
Here   for $y\in\phi(R^k)$,  $x_y\in\R^k$ is  any fixed solution of the equation $y=Ax_y+y_0$.

 When $m=n$, then $A_n=\sqrt{\mbox{det}(AA^T)}$ and the map $\phi$ is surjective i.e. $\phi(\R^k)=\R^n$. In this case theorem \ref{teo k>n} implies that $p_y\ll\mathcal L ^n$ i.e.  $Y$ is an absolutely  continuous random variable. If moreover  $k=n$ and  $A\in\R^{k\times k}$ is not-singular  then $A_k=|\mbox{det } A|$ and in this case we have 
\begin{align*}
f_Y(y)=
\frac{1}{|\mbox{det } A|}\int_{Ax+y_0=y}f_X(x)\,d\Ha^{0}(x)=\frac{f_X(A^{-1}(y-y_0))}{|\mbox{det } A|}, &\quad \text{for}\quad y\in \R^k.
\end{align*} 
\end{esem}

\begin{esem}[Sum of variables and Sample mean]
Let $\left(\Omega,\Sigma,p\right)$ be a probability measure space, let $k\in\N$ and  let $X=\left(X_1,\dots, X_k\right):\Omega\to \R^k$ be an absolutely  continuous random variable with probability density function $f_X$. We employ  Corollary \ref{teo k>1} with  
\begin{align*}
\phi:\R^k\to  \R,\quad  \phi(t)=\sum_{i=1}^k t_i,\quad  J_1(t)=|\nabla\phi|=\sqrt k.
\end{align*}
Then  the statistic $Y:=\sum_{i=1}^k X_i$ is an absolutely  continuous random variable  and its probability density function  $f_Y$  satisfies  
\begin{align*}
f_Y(y)=
\int_{\sum_{i=1}^k x_i=y}\frac{f_X(x)}{\sqrt k}\,d\Ha^{k-1}(x), &\quad \text{for a.e.}\quad y\in \R.
\end{align*} 
Let us  set $x^{k-1}:=\left(x_1,\dots,x_{k-1}\right)$ and let $\psi(x^{k-1})=\left(x^{k-1},\,y-\sum_{i=1}^{k-1}x_i\right)$ be a parametrization of the hyperplane  $\sum_{i=1}^k x_i=y$.   Using the area formula \eqref{parametrized manifold}, the last integral becomes
\begin{align*}
f_Y(y)=
\int_{\R^{k-1}}f_X\Big(x^{k-1},\,y-\sum_{i=1}^{k-1}x_i\Big)\,dx^{k-1}, &\quad \text{for a.e.}\quad y\in \R.
\end{align*}
In the particular case $k=2$ and if $X_1$, $X_2$ are independent, the last formula gives the well known convolution form for the distribution of the random variable $X_1+X_2$:
\begin{align*}
f_{X_1+X_2}(y)=
\int_{\R}f_{X_1}(t)f_{X_2}(y-t)\,dt, &\quad \text{for a.e.}\quad y\in \R,
\end{align*}
where  $f_{X_1}, f_{X_2}$ are respectively the density function of the distribution generated by $X_1,X_2$.

%
%
%If in particular $X_1,\dots,X_k$ are identically distributed and independent with common probability density function  $f:\Omega\to\R$, then 
%\begin{align*}
%f_Y(y)=
%\frac 1 {\sqrt k}\int_{\sum\limits_{i=1}^k x_i=y}\;\prod_{i=1}^k f(x_i)\,d\Ha^{k-1}(x), &\quad \text{for a.e.}\quad y\in \R.
%\end{align*}  
Moreover if  $X_1,\dots,X_k$ are identically distributed and independent with common probability density function  $f:\Omega\to\R$, then    (using also Example \ref{Affine transformations} with $\phi(x)=\frac 1 kx$), the density function of  the sample mean $Z:=\frac{1}{k}\sum_{i=1}^k X_i$ is
\begin{align*}
f_Z(y)= k\, f_Y\left(ky\right)= k
\int_{\sum\limits_{i=1}^k x_i= k y}\;\prod_{i=1}^k f(x_i)\,d\Ha^{k-1}(x), &\quad \text{for a.e.}\quad y\in \R.
\end{align*}
\end{esem}

\begin{esem}[Product and ratio of random variables]\label{ratio}
Let $\left(\Omega,\Sigma,p\right)$ be a probability measure space and  let $X:\Omega\to \R^2,\, X=\left(X_1,X_2\right)$ be an absolutely  continuous random variable with probability density function $f_X$. 
\begin{itemize}
\item[(i)] Let us  employ  Corollary \ref{teo k>1} with  
\begin{align*}
\phi:\R^2\to  \R,\quad  \phi(x_1,x_2)=x_1x_2,\quad  J_1\phi(x_1,x_2)=|\nabla\phi(x_1,x_2)|=\sqrt{x_1^2+x_2^2}.
\end{align*}
Then  the statistic $X_1X_2$ is an absolutely  continuous random variable  whose  probability density function    satisfies
\begin{align*}
f_{X_1X_2}(y)=
\int_{x_1x_2=y}\frac{f_X(x_1,x_2)}{\sqrt {x_1^2+x_2^2}}\,d\Ha^{1}(x_1,x_2)=\int_{\R\setminus\{0\}}f_X\Big(t,\frac y t\Big)\frac{1}{|t|}\,dt, &\quad \text{for a.e.}\quad y\in \R,
\end{align*} 
where we parametrized  the hyperbole $x_1x_2=y$ by $\psi(t)=\left(t,\frac y t\right)$ and we used Formula \eqref{parametrized manifold} to evaluate  the last integral.

\item[(ii)] Let us suppose $X_2\neq 0$ a.e. and let us employ Corollary \ref{teo k>1} with  
\begin{align*}
\phi:\R^2\to  \R,\quad  \phi(x_1,x_2)=\frac{x_1}{x_2},\quad  J_1\phi(x_1,x_2)=|\nabla\phi(x_1,x_2)|=\frac 1 {x_2^2}\sqrt{x_1^2+x_2^2}.
\end{align*}
Then  the statistic $\frac{X_1}{X_2}$ is an absolutely  continuous random variable  whose  probability density function   satisfies
\begin{align*}
f_{\frac{X_1}{X_2}}(y)=
\int_{\frac{x_1}{x_2}=y}f_X(x_1,x_2)\frac{x_2^2}{\sqrt {x_1^2+x_2^2}}\,d\Ha^{1}(x_1,x_2)=\int_{\R}f_X\Big(ty, t\Big)|t|\,dt, &\quad \text{for a.e.}\quad y\in \R,
\end{align*} 
where we parametrized  the line $x_1=yx_2$ by $\psi(t)=\left(ty,t\right)$ and we used  \eqref{parametrized manifold} to evaluate  the last integral.
\item[(iii)]Let $X:\Omega\to \R$ be an absolutely  continuous random variable such that $X\neq 0$ a.e. and let  $f_X$ its probability density function. We employ  Corollary \ref{teo k=n} with  \begin{align*}
\phi:\R\setminus\{0\}\to \R\setminus\{0\},\quad \phi(t)=\frac 1 t,\quad J_1(t)=|\phi'(t)|=\frac 1{t^2}.
\end{align*}
Then  the statistic $\frac 1 {X}$ is an absolutely  continuous random variable  whose probability density function satisfies    
\begin{align*}
f_{\frac 1 X}(y)=f_X\left(\frac 1 y\right)\frac{1}{y^2},\quad \text{for any}\quad y\neq 0.
\end{align*}
\end{itemize}

\end{esem}
\subsection{Order Statistics}
Let  $S_k$ be  the set of all the permutations of the set $\{1,\dots,k\}$. Let  $X=\left(X_1,\dots, X_k\right):\Omega\to \R^k$ be a  random variable and let us consider the map
$$\phi:\R^k\to\R^k,\quad x=(x_1,\dots, x_k)\mapsto (x_{(1)},\dots, x_{(k)}),$$
which associates to any vector $x$ its increasing rearrangement $(x_{(1)},\dots, x_{(k)})$  i.e. $x_{(1)}\leq \dots\leq x_{(k)}$.  The random variable $\phi\circ X:= \left(X_{(1)},\dots, X_{(k)}\right)$ is the random vector of the so-called \emph{Order Statistics} of $X$. In what follows, as an easy application of the results of the previous sections, we deduce  their  well known density functions. We start with the following Lemma which shows, in particular,  that $\phi$ has unitary Jacobian. 
 
 \begin{lem}
 Let $n\in\N$ such that $n\leq  k$, let  $I=\{i_1,i_2,\dots,i_n\}\subseteq\{1,\dots,k\}$ a subset of indexes, where $|I|=n$ and $i_1<i_2<\dots<i_n$. Let 
 $$\phi_I:\R^k\to\R^n,\quad x=(x_1,\dots, x_k)\mapsto (x_{(i_1)},\dots, x_{(i_n)})$$
 where $(x_{(i_1)},\dots, x_{(i_n)})$ is the vector obtained by extracting from the increasing rearrangement $(x_{(1)},\dots, x_{(k)})$ the component corresponding to the indexes of $I$ (note that when $I=\{1,\dots,k\}$ then  $\phi_I=\phi$).  Then the $n$-Jacobian $J_n\phi_I$ of $\phi_I$ satisfies $J_n\phi_I=1$.
 \end{lem}
 {\sc{Proof.}} Let us suppose, without any loss of generality, $I=\{1,2\dots,n\}$. For every fixed permutation $\sigma\in S_k$, let us consider the Borel subset $A_\sigma\subseteq\R^k$ defined by
 \begin{align*}
 A_\sigma:=\{x=(x_1,\dots,x_k)\in\R^k\;:\; x_{\sigma(1)}<x_{\sigma(2)}<\dots< x_{\sigma(k)}\}.
 \end{align*}
The complementary set of $\bigcup_{\sigma\in S_k} A_\sigma$ is  the set 
$$F=\{x=(x_1,\dots,x_k)\in\R^k\;:\; \exists\, i\neq j \text{ s.t. } x_i=x_j \}$$
which satisfies  $\mathcal{L}^k(F)=0$; therefore the collection of all  $A_\sigma$ is a finite  disjoint covering  of a.e point  of $\R^k$.

Let us fix $\sigma\in S_k$; on $A_\sigma$, $\phi_I$ is injective and it  coincides with the permutation of indexes operator $T_\sigma$ defined by $\sigma$:
\begin{align*}
\phi_I(x):=T_\sigma(x)=\left(x_{\sigma(1)},\dots, x_{\sigma(n)}\right),\quad \forall x\in A_\sigma.
\end{align*}
The  Jacobian matrix of $T_\sigma$ at any point $x\in A_\sigma$ is then  a  permutation  of the matrix $\left(\begin{array}{c}
  I_n  \\ 
  \hline
  0 
 \end{array}\right)$.  Recalling Definition \ref{k-Jac}, this implies in particular $J_n\phi=1$ on $A_\sigma$ which, by the arbitrariness of $\sigma$,  proves the required claim.
 \\\qed
  Let now $X_1,\dots,X_k$ be a sequence of  absolutely continuous, independent and identically distributed random variables. Let $f$ be the common  density function of each $X_i:\Omega\to	\R$ and let $F(y):=\int_{-\infty}^yf(t)\,dt$ be the associated distribution function. In the following Proposition we compute  the density of the Order statistics of the vector  $X=\left(X_1,\dots,X_k\right):\Omega\to \R^k$. We remark that $X$ is an absolutely  continuous random variable having    density function given by $f_X(y)=\prod_{i=1}^k f(y_i)$.
 \begin{prop}
\begin{itemize}
\item[(i)] The density function of the    distribution of the vector $Y=\left(X_{(1)},\dots, X_{(k)}\right)$ of the all order statistics satisfies
 \begin{align*}
f_Y(y)=k!\prod_{i=1}^kf(y_i),\quad \forall y\in \R^k \text{ s.t. } y_1<y_2<\dots<y_k
\end{align*}
and it is $0$ otherwise.
\item[(ii)] Let $i\in\N$ such that $i\leq k$.  The density function of the   distribution of the order statistic $X_{(i)}$ satisfies
\begin{align*}
f_{X_{(i)}}(y)=i {{k}\choose{i}}\,f(y)\,F(y)^{i-1}\,\Big(1-F(y)\Big)^{k-i}, \quad \forall y\in\R.
\end{align*}
\item[(iii)] Let $i,j\in\N$ such that $i<j\leq k$.  The density function of the  distribution of the vector of the two order statistics $\left(X_{(i)},X_{(j)}\right)$ satisfies
 \begin{align*}
f_{\left(X_{(i)},X_{(j)}\right)}(y_1,y_2)=&\frac{k!}{(i-1)!(j-i-1)!(k-j)!}f(y_1)f(y_2)\\[1ex]
&\times F(y_1)^{i-1}\Big(F(y_2)-F(y_1)\Big)^{j-i-1}\Big(1-F(y_2)\Big)^{k-j},\quad \forall y\in \R^2 \text{ s.t. } y_1<y_2
\end{align*}
and it is $0$ otherwise.
\end{itemize}
 \end{prop}
 {\sc{Proof.}} Using Corollary \ref{teo k=n} and the previous Lemma, we get that  the random vector $Y=\left(X_{(1)},\dots, X_{(k)})\right)$ of the Order Statistics of $X$    is an absolutely  continuous random variable   and its probability density function  $f_Y$  satisfies  
\begin{align*}
f_Y(y)&=\sum_{\sigma\in S_k}f_X\left(T_\sigma^{-1}(y)\right)=\sum_{\sigma\in S_k}f_X\left(T_{\sigma^{-1}}(y)\right)=\sum_{\sigma\in S_k}\prod_{i=1}^kf(y_{\sigma^{-1}(i)})\\
&=\sum_{\sigma\in S_k}\prod_{j=1}^kf(y_j)=k!\prod_{j=1}^kf(y_j),\quad \forall y\in \R^k \text{ s.t. } y_1<y_2<\dots<y_k.
\end{align*}
and it is $0$ otherwise. This  proves (i).

 Let $i\in\N$ such that $i\leq k$. Claim (ii)  can be proved  directly by applying, as in the previous step,   Theorem \ref{teo k>n} and the previous Lemma or alternatively by integrating the joint density $f_Y$. Indeed if we write $\hat{x_i}=\left(x_1,\dots x_{i-1}, x_{i+1},\dots, x_k\right)\in\R^{k-1}$ to denote the variable $x$ without the $x_i$ component and   if we set  $F=\{\hat{x_i}\in\R^{k-1} \;:\; x_1<x_2\dots<x_{i-1}<y<x_{i+1}<\dots<x_k\}$ then we  obtain
\begin{align*}
f_{X_{(i)}}(y)&=\int_{\{x\in\R^k:x_i=y\}}f_Y\left(x\right)\,d\hat{x_i}= k!f(y)\int_{F}\prod_{j\neq i}f(x_j)\,d\hat{x_i}\\
&= k!f(y)\int_{\{x_1<\dots<x_{i-1}<y\}}\prod_{j=1}^{i-1}f(x_j)\,dx_1\dots dx_{i-1}\int_{\{y<x_{i+1<}\dots<x_{k}\}}\prod_{j=i+1}^{k}f(x_j)\,dx_{i+1}\dots dx_{k}.
\end{align*}
Since the integrand of the first integral of the right hand side of the last equation is invariant under any  permutations of its variables then
\begin{align*}
\int_{\{x_1<\dots<x_{i-1}<y\}}\prod_{j=1}^{i-1}f(x_j)\,dx_1\dots dx_{i-1}=\frac 1{(i-1)!}\int_{\{x_j<y,\, \forall j\leq i-1\}}\prod_{j=1}^{i-1}f(x_j)\,dx_1\dots dx_{i-1}\\
=\frac 1{(i-1)!}\prod_{j=1}^{i-1}\int_{-\infty}^y f(x_j)dx_j=\frac 1{(i-1)!}F(y)^{i-1}.
\end{align*}
Analogously one has
\begin{align*}
\int_{\{y<x_{i+1<}\dots<x_{k}\}}\prod_{j=i+1}^{k}f(x_j)\,dx_{i+1}\dots d_{k}=\frac 1{(k-i)!}\left(1-F(y)\right)^{k-i}.
\end{align*}
This gives 
\begin{align*}
f_{X_{(i)}}(y)=k!\frac 1{(i-1)!(k-i)!}f(y)F(y)^{i-1}\left(1-F(y)\right)^{k-i}=i {{k}\choose{i}}f(y)F(y)^{i-1}\left(1-F(y)\right)^{k-i}
\end{align*}
which is the required claim.

The proof of (iii) follows similarly.
\\\qed

\begin{os}
	 (i) The joint density function of three or more order statistics could be derived using similar  arguments.
	
	(ii) The same methods applies also when the random variables $X_1,\dots,X_k$ are independent  but not identically distributed. For example let  $f_i(y)$ be the density function of  $X_i:\Omega\to	\R$.  Then in this case $f_X(y)=\prod_{i=1}^kf_i(y_i)$ for every $y\in\R^k$. If $Y=\left(X_{(1)},\dots, X_{(k)})\right)$ then  one obtains as before
	\begin{align*}
		f_Y(y)=\sum_{\sigma\in S_k}\prod_{i=1}^kf_i(y_{\sigma(i)}),\quad \forall y\in \R^k \text{ s.t. } y_1<y_2<\dots<y_k.
	\end{align*}
\end{os}

\subsection{Normal distributions}
Let $a\in\R$, $\sigma> 0$ and let $X:\Omega\to \R$ be a random variable. $X$ is said to have a 
\emph{ Normal (or Gaussian) distribution}
 $p_X$, and we write $p_X\sim  \mathcal N\left(a,\sigma^2\right)$, if  $p_X$ has density
 \begin{align*}
 f_X(t)=\frac{1}{\sqrt{2\pi\sigma^2}}
 \exp\left(-\frac{|t-a|^2}{2\sigma^2}\right),\quad t\in\R.
 \end{align*}
When $\sigma=0$ we also write $p_X\sim  \mathcal N\left(a,0\right)$ with the understanding  that  $p_X$ is the dirac measure $\delta_a$ at the point $a$. 
The parameters $a$ and $\sigma^2$ are called the mean and the variance of $X$, respectively.

Let $k\in\N$,  $a\in\R^k$ and let $\Sigma\in\R^{k,k}$ be a symmetric, and positive semi-definite matrix. A random variable $X=\left(X_1,\dots,X_k\right):\Omega\to \R^k$  is said to have a \emph{(multivariate) normal distribution}  $\mathcal{N}(a,\Sigma)$  if 
\begin{align*}
\langle\lambda, X\rangle \sim \mathcal{N}\Big(\langle\lambda, a\rangle, \langle\Sigma \lambda,\lambda \rangle\Big),\quad \forall \lambda\in\R^k
\end{align*}
(we write $\langle \lambda,\mu \rangle=\sum_{i}\lambda_i\mu_i$ to denote the inner product of $\R^k$).
Here $a:=E\left(X\right)$ is the mean vector and $\Sigma=\left(\sigma_{ij}\right)_{i,j}:=\mbox{Cov}(X)$ is the covariance matrix of $X$ i.e. $\sigma_{i,j}=\mbox{Cov}\left(X_i,X_j\right)$. The following very well known properties about Gaussian vectors are direct consequences of their definition  (see for example \citep[Chapter 1]{bogachev}).

\begin{prop}\label{Basic Normal}
Let  $X=\left(X_1,\dots,X_k\right):\Omega\to \R^k$ be a random variable such that $X\sim\mathcal{N}\left(a,\Sigma\right)$. 
\begin{itemize}
\item[(i)] The mean vector $a$ and the Covariance matrix $\Sigma$  uniquely characterized the Gaussian measure $p_X$.
\item[(ii)]  $X_i$, $X_j$ are independent if and only if $\sigma_{ij}=\mbox{Cov}(X_i,X_j)=0$. 
\item[(iii)] For every matrix $A\in\R^{m,k}$ one has $AX\sim\mathcal{N}\left(\langle Aa\rangle, A\Sigma A^t\right)$.
\item[(iv)]  When $\Sigma$ is positive definite  we say that $p_X$ is \emph{not-degenerate}: in this case  $X$ is  absolutely continuous and has density function 
\begin{align*}
f_X(x)=\frac{1}{(2\pi)^{\frac k 2}\mbox{det}(\Sigma)^\frac 1 2}\exp{\left(-\frac{|\Sigma^{-\frac 1 2}(x-a)|^2}{2}\right)}.
\end{align*}
\end{itemize}  
\end{prop}

\medskip
In the following Proposition we show  that, when the Covariance matrix $\Sigma$ is degenerate,  $p_X$ has a density function  with respect to the  Hausdorff measure $\Ha^m$ on  some hyperplane of $\R^k$. In what follows we say that  that a matrix $P\in \R^{k,m}$, with  $m\leq k$,  is  orthogonal if it has orthonormal columns; in this case $|Qy|=|y|$ for every $y\in\R^m$.

\begin{prop}
Let $k\in\N$,  $a\in\R^k$ and let $\Sigma\in\R^{k,k}$ be a positive semi-definite matrix with  $m=\mbox{rank} (\Sigma)\geq 1$. Let  $X=\left(X_1,\dots,X_k\right):\Omega\to \R^k$ be a random variable. Then  $X\sim\mathcal{N}(a,\Sigma)$ if and only if there exists an orthogonal matrix $P\in \R^{k,m}$ and $m$ independent random variables $Y_1,\dots,Y_m$  which satisfies $Y_i\sim \mathcal{N}\left(0,1\right)$ for every $i=1,\dots,m$ and  such that 
\begin{align*}
X=\Sigma^{\frac 1 2 }P\,Y+a,\quad Y=\left(Y_1,\dots,Y_m\right).
\end{align*}
Moreover the probability measure $p_X$ has a density function $f_X$ with respect to the  Hausdorff measure $\Ha^m$ on  the hyperplane 
$$\Sigma^{\frac 12}P\left(\R^m\right)+a=\{x\in\R^k\;:\;x=\Sigma^{\frac 12}Py+a,\;y\in\R^m\}$$
 which  satisfies   
\begin{align*}
f_X(x)=\frac{1}{(2\pi)^{\frac m 2}\Sigma^{\frac 12}_m}\exp{\left(-\frac{|y|^2}{2}\right)}, \quad x=\Sigma^{\frac 12}Py+a,
\end{align*}
where  $\Sigma^{\frac 12}_m=\prod_i{\sqrt \lambda_i}$ and  the product runs over all positive eigenvalues of $\Sigma$ (counted with their multiplicities).

\end{prop}
{\sc{Proof.}} Let us prove the first claim and let  us suppose, preliminarily, that the Covariance matrix $\Sigma$ is a diagonal matrix and, without any loss of generality, let us assume that its entries in the main diagonal are 
$$\left(\sigma_{11},\dots,\sigma_{mm}, 0,\dots,0\right),$$
 where  $\sigma_{ii}>0$ for $i\leq m$. Then from Proposition \ref{Basic Normal}, $X_1,\dots, X_m$ are independent and $X_i\sim\mathcal{N}(a_i,\sigma_{ii})$; moreover $X_i=a_i$ a.e. for $i>m$. The required claim then immediately follows  setting $Y=\left(Y_1,\dots, Y_m\right)$, with $Y_i=\frac{X_i-a_i}{\sqrt\sigma_{ii}}$, and  $P=\left(\begin{array}{c}
  I_m  \\ 
  \hline
  0 
 \end{array}\right)$, where $I_m$ is the identity matrix of $\R^{m,m}$.
% \begin{align*}
%P=\left(\begin{array}{c}
%  I_m  \\ 
%  \hline
%  0 
% \end{array}\right).
%\end{align*}

 In the general case let us diagonalize  the Covariance matrix $\Sigma$: Let $Q\in\R^{k,k}$ be an  orthogonal matrix  such that $Q\Sigma Q^t=D$, where $D$ is the diagonal matrix whose entries in the main diagonal are $\left(\lambda_1,\dots,\lambda_m, 0,\dots,0\right)$, where  the  $\lambda_i>0$ are the positive eigenvalues of $\Sigma$. From Proposition \ref{Basic Normal} the vector $Z=QX$ satisfies $Z\sim\mathcal{N}\left(Qa,D\right)$; from the previous step there exists $Y=\left(Y_1,\dots,Y_m\right)\sim \mathcal{N}\left(0,I_m\right)$ such that
 \begin{align*}
 Z=D^{\frac 1 2}\left(\begin{array}{c}
  I_m  \\ 
  \hline
  0 
 \end{array}\right)Y+Qa.
 \end{align*}
 Then since $Q\Sigma^{\frac 1 2 }Q^t=D^{\frac 1 2 }$ we get 
 \begin{align*}
 X=Q^tZ=Q^tD^{\frac 1 2}\left(\begin{array}{c}
  I_m  \\ 
  \hline
  0 
 \end{array}\right)Y+a=\Sigma^{\frac 1 2} Q^t\left(\begin{array}{c}
  I_m  \\ 
  \hline
  0 
 \end{array}\right)Y+a
 \end{align*}
 and the claim follows with $P=Q^t\left(\begin{array}{c}
  I_m  \\ 
  \hline
  0 
 \end{array}\right)$.

Finally, to prove the second claim, let us apply the first step and  Example \ref{Affine transformations} with $A=\Sigma^{\frac 1 2}P$ and $y_0=a$. Then we get that $X$ has a density function $f_X$ with respect to the  Hausdorff measure $\Ha^m$ on  the hyperplane 
$$\Sigma^{\frac 12}P\left(\R^m\right)+a=\{x\in\R^k\;:\;x=\Sigma^{\frac 12}Py+a,\;y\in\R^m\}$$
 which  satisfies   
\begin{align*}
f_X(x)&=
\frac{1}{(\Sigma^{\frac 12}P)_m}\int_{\Sigma^{\frac 12}Py+a=x}f_Y(y)\,d\Ha^{0}(y), \quad \text{for}\quad x\in \Sigma^{\frac 12}P\left(\R^m\right)+a.
\end{align*} 
Since $P$ has orthogonal columns then from Definition \ref{k-Jac} we have $(\Sigma^{\frac 12}P)_m=\prod{\sqrt \lambda_i}:=\Sigma^{\frac 12}_m$, where the product runs over all positive eigenvalues of $\Sigma$. Moreover since $\Sigma^{\frac 12}P$ has maximum  rank, the equation  $x=\Sigma^{\frac 12}Py+a$ has a unique solution. Then
\begin{align*}
f_X(x)=\frac{1}{(2\pi)^{\frac m 2}\Sigma^{\frac 12}_m}\exp{\left(-\frac{|y|^2}{2}\right)}, \quad x=\Sigma^{\frac 12}Py+a.
\end{align*}
\qed

\subsection{Chi-squared and Student's  distributions}
Let $X:\Omega\to \R^k$ be a Gaussian random vector whose  covariance matrix is the identity matrix $I_k$. If $X\sim\mathcal{N}\left(0,I_k\right)$ then the probability measure $p_{\chi^2(k)}$ induced by  $|X|^2$ is called  \emph{Chi-squared distribution with $k$-degrees of freedom} and we write $|X|^2\sim \chi^2(k)$.

 If $X$ is not-centred  i.e. $X\sim\mathcal{N}\left(\mu,I_k\right)$ for some $\mu\in\R^k\setminus\{0\}$, then the measure $p_{\chi^2(k,\lambda)}$  induced by $|X|^2$ is called  \emph{Non-central Chi-squared distribution with $k$-degrees of freedom and non-centrality parameter $\lambda=|\mu|^2>0$} and we write $|X|^2\sim \chi^2(k,\lambda)$.

\medskip
In the next Proposition we derive the density function of $|X|^2$. In what follows we consider the gamma function $\Gamma(r)=\int_0^\infty t^{r-1}e^{-r}\,dr$, $r>0$ (see e.g. \citep[page 255]{abramowitz+stegun}) and the modified Bessel function of the first kind $ I_{\nu }$  defined for $y>0$ as   
\begin{align*}
I_{\nu }(y)=(y/2)^{\nu }\sum _{j=0}^{\infty }{\frac {(y^{2}/4)^{j}}{j!\,\Gamma (\nu +j+1)}}=\frac{(y/2)^{\nu }}{\pi^{\frac 1 2}\Gamma\left(\nu+\frac 1 2\right)}\int_0^\pi e^{y\cos\theta}\left(\sin \theta\right)^{2\nu}\,d\theta,
\end{align*}
(see e.g. \citep[Section 9.6 and Formula 9.6.20, page 376]{abramowitz+stegun}).
\begin{prop}[Chi-squared  Distribution]\label{section chi squared}
Let $X:\Omega\to \R^k$ be a Gaussian random vector. If  $X\sim\mathcal{N}\left(0,I_k\right)$ then  the Chi-squared distribution $p_{\chi^2(k)}$ induced by $|X|^2$ has density function
\begin{align*}
f_{\chi^2(k)}(y)=\frac{1}{2^{\frac k 2}\Gamma\left(\frac k 2\right)}y^{\frac k 2-1}\exp\left(-\frac y 2\right),\quad \text{for any}\quad y>0.
\end{align*}
If $X\sim\mathcal{N}\left(\mu,I_k\right)$ for some $\mu\in\R^k\setminus\{0\}$ then, setting $\lambda=|\mu|^2>0$,  the Non-central Chi-squared distribution  $p_{\chi^2(k,\lambda)}$ induced by $|X|^2$ has density function

\begin{align*}
f_{\chi^2(k,\lambda)}(y)&=\frac{1}{2}\exp\left(-\frac{y+\lambda}2\right)\left(\frac y\lambda\right)^{\frac k 4-\frac 1 2}I_{\frac k 2 -1}\left(\sqrt{\lambda y}\right),\quad \text{for any}\quad y>0.
\end{align*}
\end{prop}
{\sc{Proof.}} 
Let $X\sim\mathcal{N}\left(0,I_k\right)$; using Example \ref{Square function} we have for any $y>0$
\begin{align*}
f_{\chi^2(k)}(y)&=\frac{1}{2\sqrt y}\frac{1}{(2\pi)^{\frac k 2}}\int_{|x|=\sqrt y}\exp\left(-\frac{|x|^2}2\right)\,d\Ha^{k-1}(x)
\\[1ex]&=\frac{1}{2\sqrt y}\frac{1}{(2\pi)^{\frac k 2}}\exp\left(-\frac{y}2\right)\Ha^{k-1}\left(\mathbb{S}^{k-1}\right)y^{\frac{k-1}2}=\frac{1}{2^{\frac k 2}\Gamma\left(\frac k 2\right)}y^{\frac k 2-1}\exp\left(-\frac y 2\right)
\end{align*}
which is the first claim. If  $X\sim\mathcal{N}\left(\mu,I_k\right)$ for some $\mu\in\R^k\setminus\{0\}$, then using Example \ref{Square function} again and the elementary  equality $|x-\mu|^2=|x|^2+|\mu|^2-2\langle x,\mu\rangle$ we have for any $y>0$
\begin{align*}
f_{\chi^2(k)}(y)&=\frac{1}{2\sqrt y}\int_{|x|=\sqrt y}\frac{1}{(2\pi)^{\frac k 2}}\exp\left(-\frac{|x-\mu|^2}2\right)\,d\Ha^{k-1}(x)\\[1ex]
&=\frac{1}{2\sqrt y (2\pi)^{\frac k 2}}\exp\left(-\frac{y+\lambda}2\right)\int_{|x|=\sqrt y}\exp\langle x,\mu\rangle\,d\Ha^{k-1}(x)
\\[1ex]&=\frac{y^{\frac {k}2-1}}{2(2\pi)^{\frac k 2}}\exp\left(-\frac{y+\lambda}2\right)\int_{|z|=1}\exp\langle \sqrt y z,\mu\rangle\,d\Ha^{k-1}(z).
\end{align*}
Since $\Ha^{K-1}_{\vert_{\mathbb{S}^{k-1}}}$ is rotationally invariant,  up to an  orthogonal transformation of $\R^k$ which maps $\frac{\mu}{|\mu|}$ to $e_1=\left(1,0,\dots,0\right)$, we can suppose $\frac{\mu}{|\mu|}=e_1$. Using $k$-dimensional spherical coordinates to evaluate the last integral then   we have 
\begin{align*}
\int_{|z|=1}\exp\langle \sqrt y z,\mu\rangle\,d\Ha^{k-1}(z)&=\int_{|z|=1}\exp\left(\sqrt{y\lambda}z_1\right)\,d\Ha^{k-1}(z)
\\[1ex]&=\Ha^{K-2}(\mathbb{S}^{K-2})\int_0^\pi \exp\left(\sqrt{y\lambda}\cos\theta\right)(\sin\theta)^{k-2}\,d\theta
%\\[1ex]&=\frac{2\pi^{\frac{k-1}{2}}}{\Gamma\left(\frac{k-1}2\right)}
%\frac{\pi^{\frac 1 2}\Gamma\left(\frac {k-1} 2\right)}{(\sqrt{\lambda y}/2)^{\frac k 2 -1 }}
%I_{\frac k 2 -1}\left(\sqrt{\lambda y}\right)
\\[1ex]&=(2\pi)^{\frac k 2}\left(\lambda y\right)^{-\frac k 4+\frac 1 2}I_{\frac k 2 -1}\left(\sqrt{\lambda y}\right).
\end{align*}
Combining the latter equalities gives the required last claim.\qed
\bigskip

 Finally let  $X,Y:\Omega\to \R$ be two independent random variables such that $X\sim\mathcal{N}\left(0,1\right)$ and  $Y\sim\chi^2(k)$. The probability measure $p_{T}$ induced by the random  variable $T=\frac{X}{\sqrt {Y/k}}$ is called a \emph{(Student's)  t-distribution with $k$-degrees of freedom}. In the next Proposition we use Corollary \ref{teo k>1} and example  \ref{ratio} in order to derive the density function of $p_{T}$.

\begin{prop}[Student's t-Distribution]
Let $X,Y:\Omega\to \R$ two independent random variables such that $X\sim\mathcal{N}\left(0,1\right)$ and  $Y\sim\chi^2(k)$.  The t-distribution $p_{T}$ induced by  $T=\frac{X}{\sqrt {Y/k}}$ has density function
\begin{align*}
f_T(y)=\frac{\Gamma\left(\frac{k+1}2\right)}{\sqrt{k\pi}\Gamma\left(\frac{k}2\right)}\left(1+\frac{y^2}{k}\right)^{-\frac{k+1}2},\quad \text{for any}\quad t\in\R.
\end{align*}
\end{prop}
{\sc{Proof.}}
Using Corollary \ref{teo k>1} with $\phi:\R^+\to\R^+, \phi(t)=\sqrt {t/k}$,  we have
\begin{align*}
f_{\sqrt {Y/k}}(y)&=\frac{2k^{\frac k 2}}{2^{\frac k 2}\Gamma\left(\frac k 2\right)}y^{ k -1}\exp\left(-\frac {ky^2} 2\right),\qquad \forall y>0.
\end{align*}
Then applying  example  \ref{ratio}  we get for $y\in \R$,
\begin{align*}
f_{T}(y)=\int_{0}^\infty f_X(ty)f_{\sqrt {Y/k}}(t)t\,dt
=\frac{2k^{\frac k 2}}{2^{\frac k 2}\Gamma\left(\frac k 2\right)\sqrt{2\pi}}\int_{0}^\infty \exp\left({-\frac{t^2(y^2+k)}{2}}\right)t^k\,dt,
\end{align*} 
which with the substitution $s=t^2\frac{y^2+k}{2}$ becomes
\begin{align*}
f_{T}(y)=\frac{1}{\Gamma\left(\frac k 2\right)\sqrt{k\pi}}\left(1+\frac{y^2}k\right)^{-\frac {k+1} 2 }\int_{0}^\infty e^{-s}s^{\frac {k-1}2}\,ds=\frac{\Gamma\left(\frac {k+1} 2\right)}{\Gamma\left(\frac k 2\right)\sqrt{k\pi}}\left(1+\frac{y^2}k\right)^{-\frac {k+1} 2 }.
\end{align*}
\qed

\bibliographystyle{apalike}
\bibliography{References}

\end{document}